\documentclass{amsart}
\usepackage{mathrsfs}
\usepackage{amsmath,amscd}
\usepackage{amsmath,amsthm,amssymb,float,multicol}
\usepackage{bibarts}
\usepackage{cite}
\usepackage[all,cmtip]{xy}
\usepackage{hyperref}

\usepackage[pdftex]{graphicx}
\usepackage[usenames]{color}
\usepackage[large]{subfigure}
\usepackage{xypic}

\DeclareMathOperator{\Aut}{Aut}

\DeclareMathOperator{\Chn}{Chn}
\DeclareMathOperator{\Hom}{Hom}

\DeclareMathOperator{\Pic}{Pic}

\DeclareMathOperator{\Gal}{Gal}

\DeclareMathOperator{\Sec}{Sec}
\DeclareMathOperator{\Spec}{Spec}
\DeclareMathOperator{\Isom}{Isom}
\DeclareMathOperator{\B}{B}

\DeclareMathOperator{\Won}{Won}

\DeclareMathOperator{\N}{\mathbb{N}}

\def\C{\mathbb{C}}
\def\F{\mathcal{F}}

\def\G{\mathbb{G}}

\def\K{\mathcal{K}}
\def\E{\mathcal{E}}
\def\L{\mathcal{L}}
\def\M{\overline{{M}}}
\def\N{\mathcal{N}}

\def\O{\mathcal{O}}
\def\P{\mathbb{P}}

\def\R{\mathfrak{R}}

\def\T{\mathcal{T}}
\def\U{\mathcal{U}}
\def\X{\mathcal{X}}

\def\Z{\mathbb{Z}}

\def\<{\langle}
\def\>{\rangle}
\def\({\left(}
\def\){\right)}

\newtheorem{theorem}{Theorem}[section]
\newtheorem{lemma}[theorem]{Lemma}

\newtheorem{prop}[theorem]{Proposition}
\newtheorem{cor}[theorem]{Corollary}
\newtheorem{question}[theorem]{Question}
\theoremstyle{definition}
\newtheorem{definition}[theorem]{Definition}
\newtheorem{construction}[theorem]{Construction}
\newtheorem{notation}[theorem]{Notation}

\newtheorem{hyp}[theorem]{Hypothesis}
\newtheorem*{ack}{Acknowledgment}

\theoremstyle{remark}
\newtheorem{remark}[theorem]{Remark}

\newtheorem{prin}[theorem]{Principle}
\numberwithin{equation}{section}

\begin{document}
\bibliographystyle{alpha}
\title{Homogeneous Space Fibrations over Surfaces}

\author{Yi Zhu}
\address{Pure Mathematics\\Univeristy of Waterloo\\Waterloo, ON N2L3G1\\ Canada}
\email{yi.zhu@uwaterloo.ca}
\curraddr{}


\date{}

\dedicatory{}

\begin{abstract}
By studying the theory of rational curves, we introduce a notion of rational simple connectedness for projective homogeneous spaces. As an application, we prove that over a function field of an algebraic surface over an algebraically closed field, a variety whose geometric generic fiber is a projective homogeneous space admits a rational point if and only if the elementary obstruction vanishes. 
\end{abstract}

\maketitle
\section{Introduction}

In this introduction, we work with varieties defined over an algebraically closed field $k$. By the work of Graber-Harris-Starr \cite{GHS} and de Jong-Starr \cite{dJS1}, any smooth separably rationally connected variety over a function field of a $k$-curve admits a rational point. One can ask a similar question over the function field $k(S)$, where $S$ is a surface. Under what conditions does a variety defined over $k(S)$ admit a rational point?

There are two difficulties to find rational points on varieties over $k(S)$. First, the class of separably rationally connected varieties is too large to admit rational points. By Tsen-Lang's theorem \cite{TsenLang}, any hypersurface of degree $d$ in the projective space $\P^n$ such that $d^2\le n$ over the function field $k(S)$ admits a rational point and the bound is sharp. This suggests that we should study varieties sharing the common geometric features with hypersurfaces in the above range. These varieties are examples of \emph{rationally simply connected varieties}, introduced by de Jong and Starr \cite{dJS}. Roughly speaking, they are varieties admitting lots of rational surfaces.

Secondly, there are Brauer-type obstructions to the existence of rational points. Since the Brauer group of $k(S)$ is not trivial in general, any Brauer-Severi variety corresponding to a nontrivial Brauer class has no rational point at all. On the other hand, the geometric generic fiber is a projective space. Such cohomological obstructions can be explained as a part of the \emph{elementary obstruction}, discovered by Colliot-Th{\'e}l{\`e}ne and Sansuc \cite{CTS}. The elementary obstruction vanishes if there is a rational point. 

Combining the above two observations, de Jong and Starr formulated the following principle.

\begin{prin}[de Jong-Starr \cite{dJS}]\label{principle}
	A rationally simply connected variety defined over $k(S)$ admits a rational point if the elementary obstruction vanishes.
\end{prin}

One piece of evidence for Principle \ref{principle} is de Jong-Starr's proof for the period-index theorem over $k(S)$ \cite{dJS2}. It is equivalent to prove that Principle \ref{principle} holds for Grassmannians. Later de Jong, He and Starr proved the following theorem.

\begin{theorem}[de Jong-He-Starr \cite{dJHS}] \label{split-dJHS}
	A projective homogeneous space of Picard number one over $k(S)$ admits a rational point if the elementary obstruction vanishes. 
\end{theorem}

The main ingredient of their work is to show that homogeneous spaces of Picard number one are rationally simply connected. Combining the work of Colliot-Th{\'e}l{\`e}ne, Gille, and  Parimala\cite{CTGP}, Serre's conjecture II over function fields of surfaces follows as a corollary. In 2008, Borovoi, Colliot-Th{\'e}l{\`e}ne, and Skorobogatov proved the following theorem.

\begin{theorem}[\cite{BCTS} Thm 3.8] 
	Assuming the period-index theorem and Serre's conjecture II for the function field $k(S)$ of a surface $S$, any homogeneous space of a connected linear $k(S)$-group admits a rational point if the elementary obstruction vanishes.
\end{theorem}

Borovoi-Colliot-Th{\'e}l{\`e}ne-Skorobogatov's theorem also gives evidence of Principle \ref{principle} for homogeneous spaces under group actions defined over the base field.

In this paper we formulate the rational simple connectedness for projective homogeneous varieties of higher Picard numbers. See Hypotheses \ref{H1}, \ref{H2}, \ref{H3}. These are geometric properties which can be checked after the base change to the algebraically closure. As an application, we prove that Principle \ref{principle} holds for projective homogeneous spaces with no assumptions on group actions.

\begin{theorem}\label{intromain}
	Let $X$ be a projective variety defined over a function field of a surface. Assume that the geometric generic fiber of $X$ is of the form $G/P$ for some linear algebraic group $G$ and parabolic subgroup $P$. Then
	$X$ admits a rational point if and only if the elementary obstruction vanishes.
\end{theorem} 

\begin{cor}[Starr]\label{serre2}
	Let $G$ be a quasisplit simply connected semisimple $k(S)$-group. Then every $G$-torsor admits a reduction of the structure group to the center of $G$.
\end{cor}

\begin{remark}
By the recent work of Starr and Xu \cite{starrxu}, the main techniques developed in this paper (cf., Theorem \ref{major} and Proposition \ref{twisting}, \ref{H1done} and \ref{H2done}) implies that Theorem \ref{intromain} holds over global function fields as well.
\end{remark}


\subsection{Sketch of the proof of Theorem \ref{intromain}}

Let $K$ be the function field $k(\P^1)$. Since the surface $S$ admits a pencil of curves over $\P^1$ under blowups, the function field $k(S)$ is the same as the function field $K(C)$. Finding a $k(S)$-rational point is equivalent to find a $K$-section of a fibration $\pi:X\to C$. 

Let $\pi:X\rightarrow C$ be a smooth family of projective homogeneous spaces over a curve $C$. The vanishing of the elementary obstruction is equivalent to the existence of a universal torsor $\T$ \cite{SK}. Theorem \ref{intromain} is proved by the following steps.

\begin{flushleft}
	\textbf{Step 1} There exists a sequence of ``canonically" chosen irreducible components $\{Z^e(X/C/K)\}_{e\ge e_0}$ in the moduli space of sections of $X/C$.\vspace{2mm}
	
	\textbf{Step 2} For each integer $e$, we can define an Abel map,
	$$\qquad\alpha_\T:Z^e(X/C/K)\rightarrow \{\text{the classifying stack of torsors over }C\text{ of degree }e\}$$
	by pullback of the universal torsor to get a torsor over $C$. This is a generalization of the classical Abel map to the intermediate Jacobian. The targets have coarse moduli spaces as a sequence of abelian varieties $\{A^e\}$ with lots of $K$-rational points. \vspace{2mm} 
	
	\textbf{Step 3} We then analyze the geometric properties of the Abel map and prove that the geometric genric fiber $F$ of $\alpha_\T$ is rationally connected. \vspace{2mm} 
	
	\textbf{Step 4} Applying the result of Graber-Harris-Starr \cite{GHS} on $F$, we have a section $\sigma:C\rightarrow X$ defined over $K$. 
\end{flushleft}

In section 2-5, we deal with Step 2. Here we generalize the notion of universal torsors to the relative setting, construct the Abel map and show its basic geometric properties. In Section 6, we define the sequence of components $Z^e$ as in Step 1 (Definition \ref{sequence}). 

In section 7 and 8, we prove Step 3 under Hypothesis \ref{H1}, \ref{H2}, \ref{H3}. See Theorem \ref{main0}. In section 9 and 10, we verify all hypotheses for projective homogeneous spaces which finishes Step 3.   

Section 11 is on discriminant avoidance which reduce the problem to treat with smooth family only. We conclude with the proof of Theorem \ref{intromain} and Corollary \ref{serre2} in section 12.

\begin{ack} The author would like to thank his thesis advisor Professor Jason Starr for introducing him to the subject and helpful suggestions. The author would also like to thank Professor Johan de Jong, Skip Garibaldi and Zhiyu Tian for helpful discussions. The author is grateful to the anonymous referee for his/her detailed comments and suggestions on the manuscript. \end{ack}

\section{Elementary Obstructions and Universal Torsors}
In this section, we first recall the elementary obstruction to the existence of rational points of varieties over fields, then generalize this construction to the relative case which gives an obstruction theory for the existence of sections. Throughout this section, we work with sheaves and cohomology in the fppf site. 

\subsection{Elementary Obstructions over a field}

The standard references for elementary obstructions are Colliot-Th{\'e}l{\`e}ne-Sansuc's original paper \cite{CTS} and Skorobogatov's book \cite{SK}. 

Let $K$ be a field. Let $X$ be a smooth projective $K$-variety and $\overline{X}$ be the base change of $X$ to the algebraic closure $\overline{K}$. Let $p: X\rightarrow\Spec K$ be the structure morphism. 

The relative Picard scheme $\Pic_{X/K}=R^1p_\ast \G_m$ is a fppf sheaf represented by a group variety over $K$ by \cite[$\text{n}^\circ 232$, 3.1]{FGA}. Let $S$ be the character group of $\Pic_{X/K}$, which is of multiplicative type over $K$. When $\Pic(\overline{X})$ is finitely generated, it is uniquely determined by $S$.

The set of isomorphism classes of $S$-torsors over $X$ is classified by the cohomology group $H^1(X,S)$. By \cite[ Th\'eor\`eme 1.5.1]{CTS}, there exists a long exact sequence of cohomological groups.
\begin{equation}\label{ses}\begin{CD}
0@>>>H^1(K, S)@>>>H^1(X, S)@>\chi>>\Hom_K(\Pic_{X/K},\Pic_{X/K})\\
@>\partial>>H^2(K, S)@>>>H^2(X, S)
\end{CD}\end{equation}

\begin{definition} Assume that $\Pic(\overline{X})$ is a finitely generated abelian group. An $S$-torsor $\T$ over $X$ is \emph{universal} if $\chi(\T)$ is the identity morphism on $\Pic_{X/K}$.
\end{definition}

\begin{definition}
	Let $Id$ be the identity morphism of $\Pic_{X/K}$. The class $e(X):=-\partial(Id)\in H^2(X,S)$ is called the \emph{elementary obstruction} of the variety $X$ over $K$.
\end{definition}

\begin{prop}\label{abs} Assume that $\Pic(\overline{X})$ is finitely generated. \label{torsor-ob}
	\begin{enumerate}
		\item The universal torsor exists if and only if the elementary obstruction $e(X)$ vanishes.
		\item If $X$ admits a $K$-rational point, then the universal torsor exists, or equivalently the elementary obstruction $e(X)$ vanishes.
	\end{enumerate}
\end{prop}
\proof
The first part follows from the long exact sequence (\ref{ses}). Since a $K$-rational point on $X$ gives a left inverse of the map $H^2(K,S)\rightarrow H^2(X,S)$ as in (\ref{ses}), the connecting map $\partial$ is the zero map. In particular, the elementary obstruction $e(X)$ vanishes.
\qed

\begin{theorem} [\cite{SK} Thm. 2.3.4]Let $X$ be a smooth projective $K$-variety. Assume that $\Pic(\overline{X})$ is a finitely generated abelian group. The class $e(X)\in H^2(X,S)$ coincides with the class of the following natural 2-fold extension of Galois modules.
	$$\begin{CD}
	1@>>>\G_{m,\overline{X}}@>>>\overline{K}(X)^\ast@>>>Div(\overline{X})@>>>\Pic(\overline{X})@>>>0
	\end{CD}$$\end{theorem}

\begin{remark} One may use the above theorem to give a general definition of elementary obstructions for smooth integral $K$-varieties without the assumption on the finite generation of Picard groups. However, we prefer this definition via universal torsors because we are mainly interested in the geometric aspect of the elementary obstruction. The finite generation of Picard groups holds for smooth projective rationally connected varieties.
\end{remark}

\subsection{Relative Universal Torsors}

\begin{hyp} \label{badhyp} Let $K$ be a field. Let $\pi:X\rightarrow C$ be a flat projective family of varieties over a smooth projective $K$-curve $C$. 
	Assume that the family satisfies the following conditions:
	\begin{enumerate}
		\item The geometric fibers of $\pi$ are reduced and irreducible. Hence by \cite[$\text{n}^\circ 232$, Thm 3.1]{FGA}, the relative Picard functor $\Pic_{X/C}$ is represented by a separated $C$-group scheme locally of finite type.
		\item Each closed subscheme of $\Pic_{X/C}$ which is of finite type is proper over $C$.
		\item The sheaves $R^1\pi_*\O_{X}$ and $R^2\pi_*\O_{X}$ are trivial and commute with base change.
		\item The geometric generic fiber of $\pi$ is smooth and simply connected, i.e. no finite \'etale cover.
	\end{enumerate}
\end{hyp}

\begin{remark}
	\begin{enumerate}
		\item Condition (2) as above is very restrictive. But it holds for smooth families by \cite[p. 232 Thm 3]{BLR}  and for families where the geometric fibers have isolated parafactorial singularities \cite[XI 3.1]{SGA2}.
		\item In characteristic zero,  by \cite{Kollarhigher1} Theorem 7.1, if the general fiber is rationally connected, the direct images $R^i\pi_*\O_X$ vanishes for $i>0$. The base change property holds if the geometric fibers have Du Bois singularities \cite[4.6]{dubois}. In particular, it holds for log canonical families \cite{lcdubois}.
		\item Koll\'ar proved that any smooth projective separable rationally connected variety over an algebraically closed field is simply connected \cite[Theorem 13]{Kollar-Fund}. Thus Condition (4) holds for projective families with general fibers smooth separable rationally connected.
		
	\end{enumerate}
\end{remark}

\begin{prop}\label{smoothS}
	Hypothesis \ref{badhyp} holds for the following families:
	\begin{enumerate}
		\item smooth families of projective homogeneous spaces;
		\item Lefschetz pencils of hypersurfaces in $\P^n$, where $n\ge 5$.
	\end{enumerate}
\end{prop}

\proof It suffices to check all the conditions in Hypothesis \ref{badhyp} for these families. For smooth families of projective homogeneous spaces, Condition (1) is trivial and Condition (2) holds by \cite[p. 232 Thm. 3]{BLR}. By proper and base change theorem \cite[III.12.9]{Hartshorne}, Condition (3) is implied by $h^1(X_t, \O)=h^2(X_t, \O)=0$ for every geometric fiber. When the fiber is the full flag variety, this follows from Kempf's vanishing theorem for line bundles \cite{kempfvanishing}. The general case then follows from the Leray spectral sequence. Since projective homogeneous spaces are rational, in particular, separably rationaly connected, Condition (4) follows from the remark as above.

For a Lefschetz pencil of hypersurfaces in $\P^n$, where $n\ge 5$, Condition (1) is trivial. Since the singular fibers of the pencil are local complete intersections of dimension $\ge 4 $, by \cite[XI, 3.13]{SGA2}, they have isolated parafactorial singularities. Thus Condition (2) follows. Vanishing of $h^1(X_t, \O)$ and $h^2(X_t, \O)$ gives Condition (3). Since every smooth hypersurface in $\P^n$ with dimension at least two is simply connected \cite[X, 3.10]{SGA2}, we have Condition (4).\qed




\begin{prop}\label{zero}
	Assuming Hypothesis \ref{badhyp}, the relative Picard functor $\Pic_{X/C}$ is represented by a torsion-free finitely generated isotrivial twisted constant $C$-group scheme.
\end{prop}


\proof  

By \cite[p. 231 Thm. 1 and Prop. 2]{BLR} and condition (3) of the Hypothesis, $\Pic_{X/C}$ is formally \'etale over $C$. Since $\Pic_{X/C}$ is of locally finite type over $C$, it is \'etale over $C$. Together with condition (2), each irreducible component of $\Pic_{X/C}$ is finite \'etale over $C$. 

Let $\eta$ be the generic point of $C$. The geometric generic fiber $\Pic_{X/C}(\overline{\eta})$ is isomorphic to a constant group scheme with coefficient group $\Z^r$. Indeed, the dimension of each connected  component of $\Pic_{X_{\overline{\eta}}/\overline{\eta}}$ is zero by the vanishing of $R^1\pi_*\O_{X}$. Hence $\Pic_{X_{\overline{\eta}}/\overline{\eta}}$ is the Neron-Severi group, which is finitely generated by the theorem of the base change \cite[XIII, 5.1]{SGA6}. The torsion-freeness follows from the fact that every torsion lines bundles gives an unramified cyclic cover and the simple connectedness of the geometric generic fiber.


Now we may choose a basis of constant sections of the group scheme $\Pic_{X_{\overline{\eta}}/\overline{\eta}}$, denoted by $v_1,\dots,v_r$. The section $v_1$ dominates a connected component of $\Pic_{X/C}$, say $B_1$. After taking the finite \'etale base change to $B_1$, $\Pic_{X/C}\times_C B_1$ is a $B_1$-group scheme equipped with a canonical section. We may take further finite \'etale base changes to get a $B$-group scheme with $r$ canonical sections. The sections induce a natural map $\Z^r\times_C B\rightarrow \Pic_{X/C}\times_C B$ between $B$-group schemes. The map is dominant by checking over the geometric generic fiber. Thus each connected component of $\Pic_{X/C}\times_C B$ is dominated by $B$ and finite \'etale over $B$. In particular each component is isomorphic to $B$. This implies that after taking the finite \'etale base change to $B$, $\Pic_{X/C}$ becomes a torsion-free finitely generated constant group scheme. Hence by definition, it is isotrivial.  \qed




Recall that there is an anti-equivalence between the category of finitely generated isotrivial twisted constant $C$-group schemes and the category of isotrivial finite type $C$-group schemes of multiplicative type via the following functors, c.f., \cite[X, 5.1, 5.6, 5.9]{SGA3}.
$$S\mapsto \hat{S}= \Hom_{C-gr}(S, \G_{m,C})$$
$$M\mapsto D(M)=\Hom_{C-gr}(M,\G_{m,C})$$
In particular, the category of torsion-free finitely generated twisted constant $C$-group schemes corresponds to the category of $C$-tori.

Assuming Hypothesis \ref{badhyp}, we now define a $C$-torus $S=D(\Pic_{X/C})$. There is the long exact sequence, which is a relative version of (\ref{ses}).  
\begin{equation}\label{ses1}\begin{CD}
0@>>>H^1(C, S)@>>>H^1(X, S)@>\chi>>\Hom_{C-gr}(\Pic_{X/C},\Pic_{X/C})\\
@>\partial>>H^2(C, S)@>>>H^2(X, S)
\end{CD}\end{equation}

Let $Id$ be the identity morphism of $\Pic_{X/C}$.
\begin{definition} Assuming Hypothesis \ref{badhyp}, the class $-\partial(Id)\in H^2(X,S)$ is called the \emph{elementary obstruction} for $p:X\rightarrow C$.  An $S$-torsor $\T$ over $X$ is \emph{universal} if $\chi(\T)$ is the identity morphism on $\Pic_{X/C}$.
\end{definition}

\begin{prop} Assuming Hypothesis \ref{badhyp}, we have the following:
	\begin{enumerate}
		\item the universal torsor exists if and only if the elementary obstruction vanishes;
		\item if the fibration $p:X\rightarrow C$ has a section, then the universal torsor exists, or equivalently the elementary obstruction vanishes.
	\end{enumerate}
\end{prop}
\proof The proof is the same as the absolute case in Proposition \ref{abs}.\qed



\section{Stable Sections and The Abel Map}\label{SS}
Let $X$ be a smooth proper $K$-variety and assume that there exists a universal torsor $\T$. Then there is a natural classifying map:
$$\alpha_\T:X(K)=\{K\text{-rational points on }X\}\rightarrow H^1(K, S)$$ by pulling back the 
universal torsor \cite[2.7.2]{CTS}. Thus we have a partition of rational points on $X$ indexed by elements in the Galois cohomology group $H^1(K,S)$. This map is crucial in studying the behavior of rational points in number theory, e.g., $R$-equivalent classes \cite{CTS}.

The main purpose of this section is to generalize this map in the relative setting $\pi:X\rightarrow C$ as in Situation \ref{badhyp}. In the relative setting, the classifying map is much more interesting because it carries algebraic structures. As we will see later, there is an algebraic map from the moduli space of stable sections to certain abelian varieties, which generalizes the construction in \cite[Sec. 6]{dJHS}. 

\begin{hyp}\label{keyhyp}
	Let $\pi:X\rightarrow C$ be a flat family of proper varieties over a connected smooth projective $K$-curve $C$ satisfying Hypothesis \ref{badhyp}. Let $S$ be the relative Neron-Severi torus. Assume that the universal $S$-torsor $\T$ exists over $X$. 
\end{hyp}

Let $\Sec(X/C/K)$ be the moduli functor parametrizing families of sections of $\pi:X\rightarrow C$. The functor $\Sec(X/C/K)$ is representable by a scheme which is a countable union of quasi-projective varieties by \cite[Part IV.4.c]{FGA}.  

Let $BS_{C/K}$ be the classifying stack of $S$-torsors on $C$. When $S$ is $\G_{m,C}$, the classifying stack is the Picard stack, which is an algebraic stack of finite type by \cite[Appendix 2]{artinversal}. In \cite[Chapter 4]{behrendthesis}, he proved that the classifying stack of torsors under reductive group scheme over a $K$-curve is a smooth algebraic stack locally of finite type. 

We have a natural $1$-morphism $$\alpha_\T': \Sec(S/C/K)\rightarrow BS_{C/K}$$ by pullback of the universal torsor. Namely, given a family of sections $\sigma :C\times_K T\rightarrow X$ over a $K$-scheme $T$,  $s^*\T$ gives a family of $S$-torsors over $C$.  This is called the \emph{Abel map}.

\begin{definition}
	The \emph{stack of stable sections} of the family $\pi:X\rightarrow C$, denoted by $\Sigma(X/C/K)$, is the fiber of the stabilization morphism 
	$$\pi_* : \M_{g(C)}(X)\rightarrow \M_{g(C)}(C,[C])$$
	over the identity map $Id: C\rightarrow C$.
\end{definition}

The natural $1$-morphism $\Sec(X/C/K)\rightarrow \Sigma(X/C/K)$ is represented by open immersions of schemes. Thus the proper algebraic stack $\Sigma(X/C/K)$ is a compactification of $\Sec(X/C/K)$. It is natural to ask if the Abel map can be extended to the stack of stable sections.

\begin{prop} \label{extendabel}Assuming that Hypothesis \ref{keyhyp} holds, there exists a $1$-morphism
	$$\alpha_\T:\Sigma(X/C/K)\rightarrow BS_{C/K}$$
	extending the Abel map $\alpha_\T':\Sec(X/C/K)\rightarrow BS_{C/K}$. Without ambiguity, we call the extended map $\alpha_\T$ the Abel map.
\end{prop}

\proof A family of stable sections of $\pi:X\rightarrow C$ over a $K$-scheme $T$ is equivalent to the following commutative diagram.
$$\xymatrix{
	C' \ar[d]_f \ar[r]^>>>>>>>\sigma  &
	X\times_K T  \ar@{->}[ld]^{(\pi, Id_T)}  \\
	C\times_K T }$$
The pullback of the universal torsor gives an $S$-torsor $\T$ over $C^\prime$. 

Since $S$ is a $C$-torus, there exists an \'etale morphism $g: D\rightarrow C$ which splits $S$, i.e., $S\times_C D$ is isomorphic to $\G_{m,D}^r$. Let $D^\prime$ be the fiber product $(D\times_K T) \times_{C\times_K T} C^\prime$. 

$$\begin{CD}
D^\prime @>g^\prime >> C' \\
@Vf^\prime VV @VfVV\\
D\times_K T @>g>> C\times_K T
\end{CD}$$

By descent theory, any $S$-torsor over $C^\prime$ is equivalent to a $\G_{m,D}^r$-torsor over $D'$ satisfying the descent datum. 
Let $\E$ be the pullback of $\T$ via $g^\prime$, which is a $\G_{m,D}^r$-torsor over $D'$. In particular, $\E$ is a product $\E_1\times\cdots\times\E_r$ of $\G_{m,D}$-torsors over $D'$. Let $p_1, p_2: D^\prime\times_{C^\prime} D^\prime\rightarrow D^\prime$ be the natural projections. The descent datum is given by an isomorphism 
\begin{equation}
\phi:p_1^* \E_1\times\cdots\times p_1^* \E_r \simeq p_2^*\E_1\times\cdots\times p_2^* \E_r
\end{equation}
satisfying the cocycle condition $p_{13}^*\phi=p_{23}^*\phi \circ p_{12}^*\phi$. Let $\phi_{ij}:p_1^* \E_i\rightarrow p_2^* \E_j$ be the component-wise morphism. 

Now we apply the functor $\det(Rf^\prime_*)$ to each factor of $\E$, c.f., \cite[Def. 3.11]{dJHS} and \cite{detconstruction}. We get a $\G_{m,D}^r$-torsor  $\F=\det(Rf^\prime_*\E_1)\times\cdots\times\det(Rf^\prime_*\E_r)$ over $D\times_K T$. It is easy to check that $\F$ is well-defined. 

We need to check that the torsor descends. First we construct an isomorphism $\psi: p_1^* \F \simeq p_2^* \F$. Since the functor $\det(Rf^\prime_*)$ commutes with the base change, it suffices to construct a morphism $ \psi:\det (Rf^\prime_*p_1^* \E_1)\times\cdots\times\det(Rf^\prime_*p_1^*\E_r) \rightarrow \det (Rf^\prime_*p_2^*\E_1)\times\cdots\times \det (Rf^\prime_*p_2^* \E_r)$. This can be defined component-wise by $\det(Rf^\prime_*\phi_{ij})$. Write $\psi$ as $\det(Rf^\prime_*\phi)$. To check that  $\psi$ is an isomorphism, define the inverse $\det(Rf^\prime_*\phi^{-1})$ as above and their composition is just the matrix multiplication $\det(Rf^\prime_*\phi^{-1})\circ\det(Rf^\prime_*\phi^{-1})=\det(Rf^\prime_* Id)=Id$. The descent cocycle condition follows directly from the descent cocycle condition for $\phi$ and the base change property of $\det(Rf^\prime_*)$. Therefore $\F$ descents to an $S$-torsor over $C$. 

When $C^\prime$ is $C\times T$, the construction is the same as pullback of the universal torsor, which coincides with the Abel map. \qed


\section{Rational Curves on Homogeneous Spaces}

Let $k$ be an algebraically closed field of characteristic zero. Let $X$ be a projective homogeneous space under a linear algebraic $k$-group. By Bruhat decomposition, the Picard lattice of $X$ is freely generated by the line bundles associated to the Schubert varieties of codimension one, denoted by $\L_1,\cdots,\L_r$. The effective cone is generated by $\L_i$'s. Indeed, any effective divisor $\sum_{i=1}^r a_i\L_i$ intersects each Schubert curve non-negatively by homogeneity. Thus by the intersection pairing, $a_i$'s are all non-negative. By homogeneity again, we see that the effective cone coincides with the nef cone. Thus the invertible sheaf $\L=\L_1+\dots+\L_r$ is ample. Since $X$ is simply connected and homogeneous, by Stein factorization, the invertible sheaf $\L$ is in fact very ample. We introduce some special curve classes on the projective homogeneous space $X$. 

\begin{definition}
	\begin{enumerate}
		\item The \emph{degree} of a curve $C$ in $X$ is the $\L$-degree of $C$.
		\item The degree one curves in $X$ are called \emph{lines}. 
		\item A curve (class) is \emph{simple} if $\L_i$-degree is either zero or one for all $i$'s.
		\item A curve (class) is \emph{maximal} if $\L_i$-degree is one for all $i$'s.
	\end{enumerate}
\end{definition}

Note that any stable rational curve with a simple curve class type is automorphism-free. The following result is a simple corollary of the main theorems in \cite{FP}, \cite{KimRahul}.

\begin{prop}\label{simpleisunob}
	Let $\beta$ be a simple curve class in $X$. The Kontsevich moduli space $\M_{0,n}(X,\beta)$ of pointed stable rational curves in $X$ is a fine moduli space, represented by a nonempty smooth projective rational variety.\qed
\end{prop}



\section{The Abel Sequences}

\begin{notation} \label{hyp}Let $K$ be a field of characteristic zero. Let $C$ be a smooth connected $K$-curve. Let $\pi: X\rightarrow C$ be a smooth family of projective homogeneous spaces. Assume that the relative Picard number, i.e., the rank of $\Pic_{X/C}(C)$ is one. Assume that the Picard number of the geometric generic fiber of $\pi$ is $r$. Let $S$ be the character $C$-group scheme of $\Pic_{X/C}$. Assume that the relative universal $S$-torsor $\T$ exists for the family. \end{notation}

By Proposition \ref{smoothS}, the relative Picard scheme $\Pic_{X/C}$ is a torsion-free finitely generated isotrivial twisted constant $C$-group scheme. Thus the character group scheme $S$ is an isotrivial $C$-torus. 

Let $\overline{\eta}$ be the geometric generic point over $C$. We can choose a canonical basis of the constant group scheme $\Pic_{X_{\overline{\eta}}/\overline{\eta}}$, denoted by $\L_1,\dots,\L_r$ such that $\L_i$'s are line bundles of $X_{\overline{\eta}}$ associated to the Schurbert cells of codimension one.  

By \cite[Expos\'e X Cor. 1.2 and Cor. 5.7]{SGA3}, the group scheme $\Pic_{X/C}$ is equivalent to specifying the geometric fiber at $\overline{\eta}$ as a discrete continuous $\pi_1(C,p)$-module, where $p$ is a geometric point of $C$.

\begin{lemma}\label{lemperm}
	The geometric fiber of $\Pic_{X/C}$ at $\overline{\eta}$ is  a discrete continuous permutation $\pi_1(C,\overline{\eta})$-module with the Galois invariant basis $\L_1,\dots,\L_r$.
\end{lemma}

\proof It is well known that the geometric generic fiber of $\Pic_{X/C}$ at is a discrete continuous permutation $\Gal(\overline{\eta}/\eta)$-module with the Galios invariant basis $\L_1,\cdots,\L_r$, c.f., \cite[Proof of Lemma 5.6]{CTGP}. The lemma follows from the fact that the natural map $\Gal(\overline{\eta}/\eta)\rightarrow\pi_1(C,\eta)$ is surjective by \cite[Expos\'e V Prop. 8.2]{SGA1}.\qed

\begin{construction}\label{conD}
	Since the rank of $\Pic_{X/C}(C)$ is one, $\L_1,\cdots,\L_r$ over $\overline{\eta}$ dominate a unique connected component of $\Pic_{X/C}$, denoted by $D$. By Proposition \ref{zero}, $D$ is a curve finite \'etale over $C$. Denote the structure map $D\rightarrow C$ by $\phi$. 
	
	In fact, $D$ admits the following Galois module interpretation. We choose a connected finite Galois cover $g:\widetilde{D}\rightarrow C$ which completely splits $\L_1, \cdots, \L_r$ with the Galois group $\Gamma$. In particular, $\Gamma$ acts on the set $\{\L_1,\cdots, \L_r \}$ transitively with the stabilizer group $\Gamma_0$ with respect to $\L_1$. Then $D$ is isomorphic to $\widetilde{D}/\Gamma_0$. Denote $\psi: \widetilde{D}\to D$ the quotient map. Furthermore, we have the following Cartesian diagram 
	$$\begin{CD}
	\coprod_{i=1}^r \widetilde{D}@>>>D\\
	@VVV @V\phi VV\\
	\widetilde{D}@>g>>C.
	\end{CD}$$
\end{construction}

The Neron-Severi torus $S$ in our setup is indeed \emph{quasisplit}.  

\begin{lemma}\cite[Lemma 3.2]{Brion2015} \label{perm}
	$S$ is isomorphic to $\R_{\phi}\G_{m,D}$.\qed
\end{lemma}





Now we introduce a natural $1$-morphism 
$$\R_\phi^{-1}:BS_{C/K}\rightarrow B\G_{m,D}$$
given by pulling back an $S$-torsor by $\phi$ to get a $\R_\phi \G
_{m,D}\times_C D$-torsor and then reducing the structure group to $\G
_{m,D}$ by the natural adjunction (projection). In fact, this is an equivalence of stacks and the inverse $1$-morphism is the Weil restriction functor $\R_\phi$, c.f., \cite[XXIV 8.2]{SGA3}. 

Let $\Pic_{D/K}$ be the relative Picard scheme and let $c:\B\G_{m,D}\rightarrow \Pic_{D/K}$ be the coarse moduli space map. Consider the Abel map defined in Proposition \ref{extendabel} and post-compose with $R_\phi^{-1}$ and the coarse moduli space map, we get the following.

\begin{definition}\label{ABEL}
	In Situation \ref{hyp}, the \emph{Abel map} for the family of homogeneous spaces $\pi: X\rightarrow C$ with respect to the universal torsor $\T$ is the composition,
	$$\begin{CD}
	\alpha_\T : \Sigma(X/C/K)@>>> BS_{C/K}@>\R_\phi^{-1}>> B\G_{m,D}@>c>>\Pic_{D/K}.\end{CD}$$
	Let $\Sigma^e(X/C/K)$ be the inverse image $\alpha_\T^{-1}(\Pic^e_{D/K})$. The number $e$ is called the $\T$-degree for the families of stable sections.
\end{definition}


Let $\sigma: C^\prime\rightarrow X$ be a stable section corresponding to a geometric point of $\Sigma^e(X/C/K)$. Then there exists a unique subcurve $C_0$ of $C^\prime$ such that  $\sigma$ restricting on $C_0$ is a honest section. The curve $C_0$ meets the rest of $C^\prime$ at finitely many points $p_1,\dots,p_\delta$. In fact, $\sigma$ is obtained by the honest section $\sigma_0$ attaching with $\delta$ stable rational curves $C_1,\dots, C_\delta$ at $p_1,\dots,p_\delta$, and the teeth lie in the fiber.

Let $q_{i,j}$ be the geometric points lying in the fiber of $\phi$ at $p_i$, where $j=1,\dots,r$.

\begin{prop}\label{degreeup}
	In Situation \ref{hyp}, let $\sigma: C^\prime\rightarrow X$ be a stable section corresponding to a geometric point of $\Sigma^e(X/C/K)$. Then there exists integers $e_{ij}$ such that the image under the Abel map is \begin{equation}\label{eqn}
	\alpha_\T(\sigma)=\alpha_\T(\sigma_0)\otimes \O_D(\Sigma_{i,j} e_{ij}q_{i,j}),
	\end{equation}
	and the set $\{e_{i1},\cdots, e_{ir}\}$ coincide with the set $\{\deg(\L_1|_{C_{i}}),\cdots,\deg(\L_r|_{C_{i}})\}$. 
	
	In particular, when we attach a vertical line to a section at $p_1$, the term  $\O_D(\Sigma_{j} e_{1j}q_{1,j})$ becomes $\O_D(q_{1,j})$ for some $j$. So the $\T$-degree increases by one. 
	When we attach a vertical maximal curve to a section at $p_1$, the term  $\O_D(\Sigma_{j} e_{1j}q_{1,j})$ becomes $\O_D(\Sigma_{j} q_{1,j})$. So the $\T$-degree increases by $r$.
\end{prop}

\proof Since Weil restriction functor is compatible with the base change, the statement can be checked by descent. Let $g:\widetilde{D} \rightarrow C$ be as in Construction (\ref{conD}). We have the following diagram.

$$
\xymatrix{
	\coprod \widetilde{D}' \ar[dd]_{\widetilde{\phi'}} \ar[rd]_{\widetilde{f'}}  \ar[rr] &&D' \ar[rd] \ar'[d][dd]^{\phi'}  \\
	&\coprod \widetilde{D} \ar[rr]^<<<<<<<<{\widetilde{g}} \ar[dd]_<<<<<<<{\widetilde{\phi}} && D \ar[dd]_<<<<<<<<\phi \\
	\widetilde{D}' \ar[rd]_{\widetilde{f}} \ar'[r][rr]_{g'}  && C'\ar[rd]_f \ar'[r][rr]_>>>>>>>{\sigma} && X \\
	&\widetilde{D}\ar[rr]_g &&C \ar[ur]_{\sigma_0}}
$$

Since $\widetilde{D}$ splits the Picard lattice, the pullback $g'^*\sigma^*\T$ is a $\G_m^r$-torsor. The torsor $\T$ being universal implies that $g'^*\sigma^*\T$ is isomorphic to the $\G_{m}^r$-torsor associated with $\L_1\times\cdots\times\L_r$ \cite[Prop. 8.1]{Peyre02}. By the construction of the extended Abel map as in Lemma \ref{extendabel}, $g^*\alpha_\T(\sigma)\cong\det(R\widetilde{f}_*\L_1)\times\cdots\times\det(R\widetilde{f_*}\L_r)$. Since $\L_1\times\cdots\times\L_r$ is isomorphic to $\R_{\widetilde{\phi'}}(\coprod \L_i)$, we have that
$$g^*\alpha_\T(\sigma)\cong\det(R\widetilde{f}_*\L_1)\times\cdots\times\det(R\widetilde{f_*}\L_r) \cong \R_{\widetilde{\phi}}(\coprod \det(R\widetilde{f'_*}\L_i)).$$
Thus the Abel image $\alpha_\T(\sigma)$ is given by descending the line bundle $\coprod \det(R\widetilde{f'_*}\L_i)$ to $D$. Since $\coprod \widetilde{D}$ is a disjoint union, it suffices to descend one line bundle $\det(R\widetilde{f'_*}\L_1)$ from $\psi:\widetilde{D}\rightarrow D$. 

We will show the case when the stable section $C'$ has only one rational curve $C_1$ attaching on $\sigma_0$ at $\sigma_0(p)$. The general case can be proved similarly. 

Choose a point $s\in g^{-1}(p)$ and let $F$ be the maximal vertical rational subcurve in $g'^{-1}(C_1)$ through $s$. Since $g$ is Galois over $C$, any vertical rational curve in $g'^{-1}(C_i)$ is expressed by $\gamma (C)$, for some $\gamma\in \Gamma$. By \cite[Lem. 6.7]{dJHS}, we have 
\begin{align}
\det(R\widetilde{f'_*}\L_1)&=\L_1|_{\widetilde{D}} \otimes \O_{{\widetilde{D}}}(\sum_{\gamma\in \Gamma}(\L_1.{\gamma(F}))\gamma(s))\\
&=\L_1|_{\widetilde{D}} \otimes \O_{{\widetilde{D}}}(\sum_{\gamma\in \Gamma}(\gamma^{-1}(\L_1).{F})\gamma(s))\label{eq}
\end{align}
where $(*.*)$ is the intersection pairing. By assumption, we know that $\Gamma$-orbit of $\L_1$ is the set $\{\L_1,\cdots, \L_r\}$. Descending (\ref{eq}) via $\psi:\widetilde{D}\rightarrow D$ gives the formula as in (\ref{eqn}).\qed


\begin{definition}
	In Situation \ref{hyp}, let $k$ be an algebraically closed field extension of $K$. A section of $\pi:X_k\rightarrow C_k$ is \emph{m-free} if for a general effective Cartier divisor $D$ of $C_k$ of degree $m$,
	$$H^1(C_k,\sigma^*N_{\sigma(C_k)/X_k}(-D))=0.$$
	A section is \emph{unobstructed} if it is $0$-free, and \emph{free} if it is $1$-free.
	A section is \emph{(g)-free} if it is ($2g(C_k)+1$)-free.
	
\end{definition}

\begin{definition}
	Let $X/C/K$ and $\T$ be as in Situation \ref{hyp}. Let $e_0$ be an integer. An \emph{Abel sequence} for $X/C/K$ is a sequence $(Z_e)_{e\ge e_0}$ of an irreducible component $Z_e$ of $\Sigma^e(X/C/K)$ which is geometrically irreducible and satisfies the following properties.
	\begin{enumerate}
		\item For every $e\ge e_0$, a general point of $Z_e$ parametrizes a $(g)$-free section.
		\item For every $e\ge e_0$, the Abel map restricted at $Z_e$
		$$\alpha_\T:Z_e\rightarrow \Pic^e_{D/K}$$
		is surjective and the geometric generic fiber is integral and rationally connected.
		\item For every $(g)$-free section $\overline{\sigma}:C\otimes_K \overline{K}\rightarrow X\otimes_K \overline{K}$ of $\T$-degree $e_0$, there exists an integer $\delta_0$ such that for every integer $\delta\ge\delta_0$, every stable section obtained by attaching $\delta$ lines in the fiber to $\overline{\sigma}$ lies in $Z_{e_0+\delta}$.
	\end{enumerate}
	A \emph{pseudo Abel sequence} is a sequence $(Z_e)_{e\ge e_0}$ as above where (2) is replaced by the weaker condition that the Abel map $\alpha_\T|_{Z_e}$ is surjective and the geometric generic fiber is integral.
\end{definition}

In Situation \ref{hyp}, we propose the following hypotheses. 
\begin{hyp}\label{H1}
	Let $t$ be a geometric point of $C$. Let $X_t$ be the geometric fiber over $t$. For any simple curve class $\beta$, the evaluation morphism 
	$$ev:\overline{M}_{0,1}(X_t,\beta)\rightarrow X_t$$
	is smooth surjective with integral rationally connected geometric fibers. 
\end{hyp}


\begin{hyp}\label{H2}For some integer $m$, the evaluation morphism for two-pointed chains of $m$ maximal rational curves,$$ev:\Chn_{2}(X/C,m\theta)\rightarrow X\times_C X$$has smooth integral rationally connected general fibers. \end{hyp}

\begin{hyp}\label{H3}[See in Definition \ref{deftw}] Let $\eta$ be the generic point of $C$. Let $X_{\overline{\eta}}$ be the geometric generic fiber of $\pi$. There exists a very twisting maximal scroll in $X_\eta$.
\end{hyp}

\begin{theorem}\label{major}
	In Situation \ref{hyp}, assume that Hypotheses \ref{H1}, \ref{H2}, and \ref{H3} hold. Then there exists an Abel sequence for $X/C/K$.
\end{theorem}

\proof By \cite[Lem. 4.11]{Starr}, to prove the existence of an Abel sequence, it suffices to prove when the base field $K$ is uncountable and algebraically closed. Now Theorem \ref{major} follows from Theorem \ref{main0}.\qed

\section{The Sequence of Components}

\begin{notation} \label{hyp1}Let $k$ be an uncountable algebraically closed field of characteristic zero. Let $C$ be a smooth connected $k$-curve. Let $\pi: X\rightarrow C$ be a smooth family of projective homogeneous spaces. Assume that the relative Picard number, i.e., the rank of $\Pic_{X/C}(C)$ is one and assume that the Picard number of each geometric fiber is $r$.
	Let $S$ be the character $C$-group scheme of $\Pic_{X/C}$. Let $\phi: D\rightarrow C$ be a finite \'etale morphism such that $S=R_\phi \G_{m,D}$ as in (\ref{conD}). Assume that the universal $S$-torsor $\T$ exists for the family. \end{notation}

\begin{lemma}[\cite{GHS}]
	Let $X/C/k$ be as in Notation \ref{hyp1}. Then there exist $(g)$-free sections. \qed
\end{lemma}

de Jong, He and Starr \cite{dJHS} introduced an important class of stable sections, the \emph{porcupines}. They are unobstructed and have nice inductive structures.
\begin{definition}
	A \emph{porcupine} in $X/C/k$ is a stable section $\sigma:C^\prime\rightarrow X$ such that
	\begin{enumerate} 
		\item the associated section $\sigma_0:C\rightarrow X$ is $(g)$-free,
		\item each vertical curve $\sigma|_{C_i}:C_i\rightarrow X_{t_i}$ is a line in the fiber of $\pi$,
		\item the attaching points of vertical curves are all distinct on $C$.
	\end{enumerate}
	We will call the section $\sigma_0$ the \emph{body}, and the vertical curves the \emph{quills}.
\end{definition}

Recall the following standard deformation results in \cite[Prop. 5.2]{Starr}.
\begin{lemma}\label{defoporc}
	\begin{enumerate}
		\item The parameter space $Porc^e(X/C/k)$ of porcupines of $\T$-degree $e$ is represented by an open smooth subscheme of $\Sigma^e(X/C/k)$.
		\item The closed subscheme $Porc^{e,\ge 1}(X/C/k)$ of $Porc^e(X/C/k)$ parametrizing porcupines with at least $1$ quill is a simple normal crossing divisor.
		\item The open subscheme $Porc^{e,\delta}(X/C/k)$ of $Porc^e(X/C/k)$ parameterizing porcupines with exactly $\delta$ quills is a smooth, locally closed subscheme of $Porc^e(X/C/k)$ of pure codimension $\delta$.\qed
	\end{enumerate}
\end{lemma}

There is a natural morphism
$$\Phi_{body}:Porc^{e,\delta}(X/C/k)\rightarrow Porc^{e-\delta,0}(X/C/k)$$
which forgets all the $\delta$ quills. Let $D_\delta$ be the $\delta$-fold symmetric product of $D$ and let $D_\delta^{\circ}$ be the dense open subset of $D_\delta$ parametrizing reduced divisors with reduced images on $C$. By Proposition \ref{degreeup}, define the refined body morphism,
$$\Phi_{body}^\prime:Porc^{e,\delta}(X/C/k)\rightarrow Porc^{e-\delta,0}(X/C/k)\times D_\delta^\circ$$
which sends a porcupine $\sigma:C^\prime\rightarrow X$ with $\delta$ quills to its body together with the attaching divisor $B_\sigma=\O_D(t_1+\dots+t_\delta)$ on $D$.

\begin{lemma}\label{body}
	In Situation \ref{hyp1}, assume that Hypothesis \ref{H1} holds. The refined body morphism
	$$\Phi_{body}^\prime:Porc^{e,\delta}(X/C/k)\rightarrow Porc^{e-\delta,0}(X/C/k)\times D_\delta^\circ$$ 
	is smooth surjective with irreducible rationally connected geometric fibers. 
\end{lemma}

\proof Given a section $\sigma$ in $Porc^{e-\delta,0}(X/C/k)$ and a reduced divisor $B=t_1+\dots+t_\delta$ in $D_\delta^\circ$, let $F$ be the space of  porcupines having the body $\sigma$ and $\delta$ quills with the attaching divisor $B$. For each $t_i$, there is a unique line class $l_i$ such that the attachment divisor is $t_i$. 
Let $F_i$ be the fiber of the evaluation morphism $\M_{0,1}(X/C,l_i)\rightarrow X$ over the point $\sigma(\phi(t_i))$.
By Hypothesis \ref{H1}, $F_i$ is a smooth integral rationally connected variety. Therefore, $F$ is the product of all $F_i$'s, which is again a smooth integral rationally connected variety. \qed 

\begin{lemma}\label{irredcomp}
	In Situation \ref{hyp1}, assume that Hypothesis \ref{H1} holds. Let $Z_{e_0}$ be an irreducible component of $\Sigma^{e_0}(X/C/k)$ whose general points parametrize $(g)$-free sections. For every $e\ge e_0$, there exists a unique irreducible component $Z_e$ such that every porcupine with body in $Z_{e_0}$ and with $e-e_0$ quills lies in $Z_e$.
\end{lemma}

\proof Let $Porc^{e_0,0}(X/C/k)_Z$ be the open subscheme of $Z_{e_0}$ parametrizing free sections. The space of porcupines with the body in $Porc^{e_0,0}(X/C/k)_Z$ and $e-e_0$ quills is irreducible by Lemma \ref{body} and unobstructed by Lemma \ref{defoporc}. Thus it is contained in a unique irreducible component of $\Sigma^{e}(X/C/k)$.\qed 

\begin{definition}\label{sequence}
	For every integer $e\ge e_0$, $Z_e$ is the \emph{distinguished irreducible component}  of $\Sigma^e(X/C/k)$ associated to $Z_{e_0}$.
\end{definition}

Combining Lemma \ref{irredcomp} and the proof of Lemma 5.7 and 5.8 in \cite{Starr}, we have the irreducibility of the geometric generic fiber of the Abel map.

\begin{prop}\label{irreducible}
	In Situation \ref{hyp1}, assume that Hypothesis \ref{H1} holds. For every $e\ge e_0+2g(D)-1$, the Abel map
	$$\alpha_\T|_{Z_e}:Z_e\rightarrow \Pic^e_{D/K}$$
	is dominant with irreducible geometric generic fiber.\qed
\end{prop}

\section{Pencils of Simple Combs}
In this section, let $X/C/k$ and $\T$ be as in Notation \ref{hyp1}.

\begin{definition} Let $\sigma$ be a free section of $X/C/k$.
	A \emph{simple $\sigma$-comb} is a stable section of $\pi:X\rightarrow C$ with the body $\sigma$ such that the vertical curves are simple stable rational curves in the fiber with distinct attaching points on $C$. 
	
	A \emph{maximal} comb is a simple comb with all the vertical curves maximal.
\end{definition}

\begin{definition}
	A two-pointed chain of rational curves in $\Sigma^e(X/C/k)$ is \emph{useful} if the marked points and the nodes parametrize unobstructed non-stacky points in $\Sigma^e(X/C/k)$. We say that the two marked points are \emph{rationally equivalent}.
\end{definition}

\begin{lemma}\label{ICsmooth}
	Any simple comb of $\T$-degree $e$ lies in the unobstructed non-stacky locus of $\Sigma^e(X/C/k)$.  
\end{lemma}

\proof For any simple comb, the body is a free section and vertical curves are free. By \cite[II.7.5]{Kollar}, the comb is unobstructed. By Proposition \ref{simpleisunob}, any vertical curves of a simple comb is non-stacky. Thus the comb itself is non-stacky. \qed

\begin{lemma}\label{movingquill}
	In Situation \ref{hyp1}, assume that Hypothesis \ref{H1} holds. Let $P\in\Sigma^e(X/C)$ be a porcupine with the body $\sigma$ and $\delta$-quills. Let $Q$ be a simple $\sigma$-comb. If the Abel images $\alpha_\T(P)$ and $\alpha_\T(Q)$ are the same, $P$ and $Q$ are rationally equivalent in $\Sigma^e(X/C)$.
\end{lemma}

\proof Since $P$ and $Q$ share the same body, by Proposition \ref{degreeup}, the attaching divisors $B_P$ and $B_Q$ are linearly equivalent divisors on $D$. Thus there exists a pencil $\P^1\rightarrow D_\delta$ connecting them. The pencil gives a rational curve in $Porc^{e-\delta,0}(X/C/k)\times D_\delta$ by the following composition.
$$\begin{CD}
\P^1 @>>> D_\delta @>(s,Id)>> Porc^{e-\delta,0}(X/C/k)\times D_\delta
\end{CD}$$
Since the attaching divisor $B_P$ is in $D_\delta^\circ$, the rational curve intersects the image of the refined body morphism
$\Phi_{body}^\prime:Porc^{e,\delta}(X/C/k)\rightarrow Porc^{e-\delta,0}(X/C/k)\times D_\delta$ by Lemma \ref{body}. 
By the result of Graber-Harris-Starr \cite{GHS}, we can lift to a rational curve in $\Sigma^e(X/C/k)$ whose general points parameterize porcupines. Specializing the family of porcupines over $B_Q$, we get a simple $\sigma$-comb $Q^\prime$ with the attaching divisor $B_Q$. Lemma \ref{ICsmooth} implies that $P$ and $Q^\prime$ are rationally equivalent. By Hypothesis \ref{H1}, $Q$ and $Q^\prime$ are connected by a useful chain of rational curves in $\Sigma^e(X/C/k)$. Therefore $P$ and $Q$ are rationally equivalent. \qed

\begin{definition}
	A \emph{maximal scroll} $R$ in $X/C$ is a morphism $r:R\rightarrow X$ such that $R\rightarrow C$ is a smooth geometrically generic ruled surface and each fiber maps to a maximal curve with at most two irreducible components. We say that $r(R)$ is the image of the maximal scroll $R$.

	A chain of $m$ maximal scrolls is \emph{transversal} if each fiber maps to a chain of $m$ maximal curves with at most $m+1$ irreducible components.
\end{definition}


\begin{lemma}\label{movingbody}
	In Situation \ref{hyp1}, assume that Hypothesis \ref{H1} holds. Assume that there exists two sections $s_0$ and $s_\infty$ on a maximal scroll $R\rightarrow C$ such that the corresponding sections $\sigma_0:=r(s_0)$ and $\sigma_\infty:=r(s_\infty)$ on $X$ are free over $C$.  Then there exists an integer $N$ such that a general maximal $\sigma_0$-comb $C$ with $N$-teeth is rationally equivalent to a simple $\sigma_\infty$-comb.
\end{lemma}

\proof For any effective divisor $D$ on $C$, let $R_D$ be the pullback divisor on $R$. When $D$ is general, $R_D$ is a disjoint union of smooth maximal curves. There exists an integer $N$ such that for a general divisor $D$ of degree $N$, the linear system $|s_0(C)+R_D|$ is sufficiently ample and the codimension one points of the linear system parametrize nodal curves, cf. \cite[Lem. 9.5]{dJHS}. In particular, the divisor $s_0(C)+R_D$ is linearly equivalent to some divisor $s_\infty(C)+E$. Since the maximal scroll contains singular fibers like a union of two simple curves, here $r(E)$ is a disjoint union of simple rational curves. Let $P$ be the maximal $\sigma_0$-comb associated to  $r(s_0(C)+R_D)$ and let $Q$ be the simple $\sigma_\infty$-comb associated to $r(s_\infty(C)+E)$.
There is a union of two general pencils joining $P$ and $Q$ such that general points parametrize nodal divisors, i.e., $P$ is rationally equivalent to $Q$. This proves Lemma \ref{movingbody} when the maximal $\sigma_0$-comb is contained in the image of $R$. For the general case, there exists a useful chain of rational curves parametrizing the family of maximal $\sigma_0$-combs by pushing all vertical maximal curves into the scroll $R$ by Hypothesis \ref{H1}.\qed

\begin{prop}\label{scrollchain}
	In Situation \ref{hyp1}, assume that Hypothesis \ref{H1} and \ref{H2} hold. Let $\sigma_0$, $\sigma_\infty$ be two $(g)$-free sections of $\pi: X\rightarrow C$. Let $T_0$, resp. $T_\infty$ be the unique irreducible component of $\Sigma(X/C/k)$ containing $\sigma_0$, resp. $\sigma_\infty$ as a smooth point. Then there exists an irreducible open subset $T\subset \Sec(\Chn_2(X/C, m\theta)/C)$ satisfying the following:
	\begin{enumerate}
		\item $T$ parametrizes a family of transversal chains of $m$ maximal scrolls;
		\item $ev_{0,\infty}|_T: T\rightarrow \Sec(X/C)\times\Sec(X/C)$ dominates $T_0\times T_\infty$;
		\item For each $\tau$ in $T$, $ev_i\circ\tau: C\rightarrow X$ gives a free section for $i=1,\dots,m-1$, where $ev_i$ is the evaluation morphism of a node on a chain.
	\end{enumerate}
\end{prop}
$$
\xymatrix{
	\Chn_2(X/C,m\theta) \ar[dd]_{Chn_2} \ar[rd]_{ev_{0,\infty}} \ar[rrd]^{ev_i} \\
	&X\times X \ar[ld]_{(\pi,\pi)} &X \ar[lld]^\pi\\
	C}
$$

\proof Consider the following commutative diagram. 
$$
\xymatrix{
	&V=\Chn_2(X/C,m\theta) \ar[ldd]_{Chn_2} \ar[d]^{ev_{0,\infty}} & C\times T \ar@{.>}[d] \ar@{.>}[l]\\
	&X\times X \ar[ld]|{(\pi,\pi)} &C\times T_0\times T_\infty \ar[lld] \ar[l]\\
	C}
$$
By \cite[Lem. 4.12, 4.17, Prop. 4.15]{dJHS}, there exists a variety $T$ parametrizing free sections of $Chn_2:V\rightarrow C$ and a dominant morphism $T\rightarrow T_0\times T_\infty$, such that the above diagram commutes.

Since $T$ parametrizes free sections and $ev_i:\Chn_2(X/C,m\theta)\rightarrow X$ is smooth, (3) follows from \cite[Lem. 3.6]{HS} Lemma 3.6. 

Finally, it suffices to show that a general section $\tau: C\rightarrow \Chn_2(X/C,m\theta)$ in $T$ gives a transversal chain of $m$ maximal scrolls. There exists a simple normal crossing divisor $\Delta$ in $\Chn_2(X/C,m\theta)$ parameterizing chains of $m$ maximal curves with at least $m+1$ irreducible components. Since $\tau$ is free, a general deformation of $\tau$ intersects the boundary strata $\Delta$ transversally by \cite[II.3.7]{Kollar}. 
\qed

\begin{prop}\label{boundaryRC}
	In Situation \ref{hyp1}, assume that Hypothesis \ref{H1} and \ref{H2} hold. Let $T_0$, resp, $T_\infty$ be an irreducible component of $\Sigma(X/C/k)$ whose general point parameterizes a $(g)$-free section of $\T$-degree $e_0$, resp, $e_\infty$. Let $Porc^{e}(X/C/k)_{T_0}$, resp, $Porc^{e}(X/C/k)_{T_\infty}$ be the moduli space of porcupines with bodies in $T_0$, resp, $T_\infty$.
	
	Then there exists an integer $E$ such that for any integer $e\ge E$ there exists a dense open subscheme $$U\subset Porc^{e}(X/C/k)_{T_0}\times_{\alpha_\T,\Pic^e_{D/k},\alpha_\T} Porc^{e}(X/C/k)_{T_\infty}$$
	in which any pair of porcupines $(P_0,P_\infty)$ are rationally equivalent in  $\Sigma^e(X/C/k)$.
\end{prop}

\proof For a general pair of $(g)$-free sections $(\sigma_0,\sigma_\infty)$, by Proposition \ref{scrollchain}, there is a transversal chain of $m$ maximal scrolls connecting them. Let $R_1,\dots,R_m$ be the maximal scrolls and let $\sigma_1,\dots,\sigma_{m-1}$ be the intermediate sections. Let $N_i$ be the integer as in Lemma \ref{movingbody} for the pair  $(R_i, \sigma_{i-1},\sigma_i)$. Choose $E= \max\{e_0,e_\infty\}+2g(D)+r\sum_{i=1}^m N_i$. For any integer $e\ge E$, let $P_0$ be a general porcupine of $\T$-degree $e$ with the body $\sigma_0$. By Lemma \ref{movingquill} and Proposition \ref{body}, $P_0$ is rationally equivalent to a general simple $\sigma_0$-comb $Q_0$ such that the teeth are the union of $N_1+\dots+N_m$ general maximal curves and lines. By Lemma \ref{movingbody}, there exists a useful chain connecting the sub-$\sigma_0$-comb of $Q_0$ with the teeth $N_1$-maximal curves and a simple $\sigma_1$-comb. The remaining teeth of $Q_0$ deform along the rational chain by Hypothesis \ref{H1}. Therefore $P_0$ is rationally equivalent to a simple $\sigma_1$-comb $P'_1$ with at least $N_2+\dots+N_m$ maximal curves. We can continue by applying Lemma \ref{movingbody} until we get a simple $\sigma_\infty$-comb $P'_\infty$. By Lemma \ref{movingquill} again, $P'_\infty$ is rationally equivalent to a general procupine $P_\infty$ having the body $\sigma_\infty$ and the same Abel image as $P_0$.\qed

\begin{cor}\label{pseudoabel}
	In Situation \ref{hyp1}, assume that Hypothesis \ref{H1} and \ref{H2} hold. Let $(Z_e)_{e\ge e_0}$ be the sequence of irreducible components of $\Sigma(X/C/k)$ defined in (\ref{sequence}). Then $(Z_e)_{e\ge e_0}$ is a pseudo Abel sequence for $X/C/k$.
\end{cor}

\proof By Lemma \ref{irredcomp} and Proposition \ref{irreducible}, it suffices to show that the sequence satisfies condition (3) of the pseudo Abel sequence. Let $\sigma$ be a $(g)$-free section. By Proposition \ref{boundaryRC}, the porcupine obtained by attaching sufficiently many quills is rationally equivalent to a porcupine in $Z_e$. Since useful chains does not leave $Z_e$, it lies in $Z_e$.   \qed

\section{Twisting Maximal Scrolls and the Abel Sequence}
In this section, let $X/C/k$ and $\T$ be as in Notation \ref{hyp1}. Let $\xi: C\rightarrow \M_{0,1}(X/C,\theta)$ be a $1$-morphism. This is equivalent to a family of pointed rational maximal curves over $C$ as the following.
$$\xymatrix{
	R \ar[d]^{p} \ar[r]^{ev}  &
	X  \ar@{->}[ld]^\pi  \\
	C \ar@/^1pc/[u]^\sigma}$$
Let $D$ be the divisor $\sigma(C)$ in $R$.

\begin{definition}
	We say that a section $s$ of $X/C$ is \emph{penned} in a maximal scroll $R$ if it coincides with the section $ev\circ\sigma(C)$ in the scroll $R$.
\end{definition}

\begin{definition}\label{deftw}
	The $1$-morphism $\xi:C\rightarrow \M_{0,1}(X/C,\theta)$ is a \emph{$m$-twisting maximal scroll} if the pair $(R,D)$ determined by $\xi$ satisfies the following properties:
	\begin{enumerate}
		\item $R$ is a maximal scroll in $X$;
		\item The sheaf $\O_R(D)$ is globally generated and non-special;
		\item The normal bundle $N_{R/X}$ is globally generated and non-special;
		\item For every divisor $\Gamma$ on $C$ of degree $\le m$, $H^1(R,N_{R/X}(-D)\otimes_{\O_R} p^*\O_C(-\Gamma)))=0$. 
	\end{enumerate}
	When $m=2$, we say that $\xi$ is \emph{very twisting maximal scroll}.
\end{definition}

\begin{prop}[\cite{Starr} Lemma 7.3]\label{surfacetocurve}
	The $1$-morphism $\xi:C\rightarrow \M_{0,1}(X/C,\theta)$ is a $m$-twisting maximal scroll if and only if it satisfies the following: 
	\begin{enumerate}
		\item $\xi(C)$ intersects the boundary divisor of $\M_{0,1}(X/C,\theta)$ transversally;
		\item The sheaf $p_*\O_R(D)$ is globally generated and non-special;
		\item The composition $ev\circ\xi:C\rightarrow X$ is a free section;
		\item The sheaf $\xi^{\ast}T_{ev}\otimes_{\O_C} \O_C(-\Gamma)$ is globally generated and non-special for every divisor $\Gamma$ on $C$ of degree $\le m$.
		
	\end{enumerate}
	When $g(C)=0$, condition (2) is equivalent to that $\xi^*T_\Phi$ is globally generated and non-special.\qed
\end{prop}

\begin{definition}
	Let $Y$ be a projective homogeneous space over algebraically closed field of characteristic zero. A maximal scroll $\zeta:\P^1\rightarrow \M_{0,1}(Y,\theta)$ is \emph{very twisting} if the induced morphism $\P^1\rightarrow \M_{0,1}(Y\times \P^1/\P^1,\theta)$ is very twisting. 
	
	A very twisting maximal scroll in $Y$ is \emph{wonderful} if both sheaves $p_*\O_R(D)$ and $p_* N_{R/X\times \P^1}$ are ample.
\end{definition}

\begin{lemma}[Lemma 12.8 in \cite{dJHS}] \label{wonderful}
	Let $Y$ be a projective homogeneous space over algebraically closed field of characteristic zero. If $Y$ has a very twisting maximal scroll, then there exist wonderful $m$-twisting maximal scrolls for arbitrary $m\ge 0$. \qed
\end{lemma}

\begin{lemma}\label{pen}
	In Situation \ref{hyp1}, assume that Hypothesis \ref{H1} holds. Every section is penned in a maximal scroll in $X/C$. 
\end{lemma}

\proof Let $\sigma$ be a section of $\pi: X\rightarrow C$. Consider the following fiber product.
$$\xymatrix{
	F \ar[r]^{} \ar[d]_{{ev'} } & \M_{0,1}(X/C,\theta) \ar[d]_{ev} \\
	C \ar[r]^\sigma &X}$$

By hypothesis \ref{H1}, $F$ is smooth over $C$ with rationally connected geometric fibers. By \cite{GHS}, there exists a section $\xi:C\rightarrow\M_{0,1}(X/C,\theta) $. By attaching sufficiently many very free curves in the fiber of $ev'$ on $\xi$, a general deformation of the comb parametrizes a free section and thus intersects the boundary strata $\Delta$ transversally by \cite[II.3.7]{Kollar}. \qed

\begin{prop}\label{twistingpen}
	In Situation \ref{hyp1}, assume that Hypothesis \ref{H1}, \ref{H2} and \ref{H3} hold. Let $(Z_e)_{e\ge e_0}$ be the pseudo Abel sequence in Corollary \ref{pseudoabel}. For every $e\ge e_0\gg0$, the irreducible component $Z_e$ contains a section $\sigma$ which is penned in a very twisting maximal scroll.
\end{prop}

\proof Let $\sigma$ be a free section in $Z_{e_0}$. By Lemma \ref{pen}, $\sigma$ is penned in a maximal scroll $R$ in $X/C$ which corresponds to a $1$-morphism $\rho:C\rightarrow \M_{0,1}(X/C,\theta) $. Deforming $\rho$ a little bit, we may assume that a general pointed rulings of $R_t$ is contained in the dense open subset of $\M_{0,1}(X/C,\theta)$ swept out by a fixed wonderful very twisting maximal scroll $g$ in some fiber of $\pi$, cf., \cite[Lem. 12.9]{dJHS}. 

Now there are arbitrarily many wonderful very twisting scrolls $g_{t_i}:\P^1\rightarrow \M_{0,1}(X_{t_i},\theta)$ such that $g_{t_i}(0)=\rho(t_i)$ and they are algebraically equivalent to $g$. Gluing $g_{t_i}$'s on $\rho$ at $\rho(t_i)$'s, we construct a comb $C\cup \cup_i g_{t_i}\rightarrow \M_{0,1}(X/C,\theta)$. By \cite[Lem. 12.11]{dJHS} and the standard comb smoothing argument, there exists $r_0$, for any $t\ge t_0$, after attaching $r$ wonderful very twisting scrolls, a general point smoothing $\xi$ of the comb corresponds to a very twisting maximal scroll in $X/C$. If the $\T$-degree of the section $\sigma_g$ in the wonderful scroll $g$ is $d$, the section $\sigma_\xi$ in the maximal scroll $\xi$ is of $\T$-degree $e_0+td$. Since the sections in $g_{t_i}$'s are free rational curves in $X_{t_i}$, the section $\sigma_\xi$ lies in $Z_{e_0+rd}$. This proves the proposition when $e=e_0+rd$.

The general case follows by repeating the above argument for sections in $Z_{e_0+1}$,$\dots$, $Z_{e_0+d-1}$. \qed 

\begin{cor} \label{perfectpen}Notations and assumptions are as in Proposition \ref{twistingpen}. Let $C_{e+r,\theta}$ be the moduli space of maximal combs with exactly one tooth and with the bodies in $Z_e$. Then a general maximal comb in $C_
	{e+r,\theta}$ is contained in the image of a very twisting maximal scroll for $e\gg 0$.
\end{cor}

\proof By Proposition \ref{twistingpen}, choose $e\gg 0$ such that a general point of $Z_e$ is contained in a very twisting maximal scroll. It suffices to show that a deformation of combs in $C_{e+r,\theta}$ can be followed by a deformation of twisting maximal scrolls containing the combs. This follows from $H^1(R,N_{R/X}(-\sigma-R_q))=0$.\qed

\begin{theorem}\label{main0}
	In Situation \ref{hyp1}, assume that Hypothesis \ref{H1}, \ref{H2} and \ref{H3} hold. For $e_0\gg 0$, the pseudo Abel sequence in Corollary \ref{pseudoabel} is an Abel sequence for $X/C/k$.
\end{theorem}

\proof By Corollary \ref{pseudoabel}, it suffices to show that for any $e\ge e_0\gg0$, the extended Abel map
$$\alpha:Z_e\rightarrow \Pic_{D/k}$$
has rationally connected geometric generic fibers. Since the target is an abelian variety, we do not worry about the rationally equivalent classes leaving the fiber of the Abel map.

We choose an integer $e_0$ such that for any $e\ge e_0$, Corollary \ref{perfectpen} holds. For any $e\ge e_0+r$, there exists an open $U_{e,\theta}\subset C_{e,\theta}$ such that every comb is contained in a very twisting maximal scroll. By \cite[Lem. 12.5]{dJHS}, every comb in $U_{e,\theta}$ is rationally equivalent to a point in the interior of $Z_e$. Since $U_{e,\theta}$ is of codimension one in $Z_e$, a general point of $Z_e$ is rationally equivalent to a general point of $U_{e,\theta}$.

Similarly, if $e\ge e_0+2r$, a general point $Z_{e-r}$ is rationally equivalent to a general point in $C_{e-r,\theta}$. Also note that the forgetting-tooth map $C_{e,\theta}\rightarrow Z_{e-r}\times C$ has rationally connected geometric fibers by Hypothesis \ref{H1}. Thus a general point in $C_{e,\theta}$ is rationally equivalent to a general point in $C_{e,2\theta}$, i.e. a general maximal comb with exactly two quills.

For any $i=0,\dots,r-1$ and for any $d\ge0$, let $e=e_0+i+dr$. By repeating the argument above, a general point in $Z_e$ is rationally equivalent to a general point in $C_{e,d\theta}$ with body in $Z_{e_0+i}$.  

By the proof of Proposition \ref{boundaryRC}, for each $i$, there exists $E_i$ such that two general points in $C_{e,d\theta}$ with the same Abel images are rationally equivalent if $d>E_i$. 

Let $E=\max_i\{E_i\}$. For any $e>e_0+rE$, given two general points in $Z_e$ with the same Abel images, each of them is rationally equivalent to a general point in $C_{e,d\theta}$. From previous paragraph, they are rationally equivalent in $Z_e$.\qed

\section{Very Twisting Maximal Scrolls on Homogeneous Spaces} 
Let $X$ be a projective homogeneous space over an algebraically closed field $k$ of characteristic zero. Let $\theta$ be the maximal curve class. Let $\zeta: \P^1\rightarrow \M_{0,1}(X, \theta)$ be a $1$-morphism. We have the following diagram,

$$\xymatrix{
	\P^1 \ar[r]^>>>>>{\zeta}  &
	\M_{0,1}(X, \theta) \ar[d]^{\Phi} \ar[r]^<<<<{ev}&
	X \\
	&\M_{0,0}(X, \theta)}$$
where $\Phi$ is the forgetful map and $ev$ is the evaluation map. By homogeneity and generic smoothness, the evaluation map $ev$ is a smooth morphism. In particular, the relative tangent bundle $T_{ev}$ is locally free. 

\begin{definition}
	The $1$-morphism $\zeta:\P^1\rightarrow \M_{0,1}(X,\theta)$ is \emph{very twisting} if the following conditions hold:
	\begin{enumerate}
		\item the vector bundle $\zeta^*T_{ev}$ is ample;
		\item the vector bundle $(ev\circ \zeta)^*TX$ is globally generated;
		\item the image $\zeta(\P^1)$ is in the smooth locus of the forgetful map $\Phi$ and the line bundle $\zeta^*T_\Phi$ is globally generated.
	\end{enumerate} 
	In this case, we say that $X$ admits \emph{a very twisting maximal scroll}. 
\end{definition}

\begin{remark}
	The definition of a very twisting $1$-morphism over any variety is given in \cite[4.3]{HS}. It is still open how to find a very twisting $1$-morphism on varieties in general. The only known examples are general low degree complete intersections in $\P^n$ and projective homogeneous spaces of Picard number one cf. \cite{dJHS}. In these cases, one can construct a very twisting scroll of the minimal curve class type. On the other hand, for varieties with higher Picard numbers, a very twisting morphism usually does not exist for minimal curve classes. Thus the existence result depends on the choice of a ``good" curve class. For smooth quadric surfaces in $\P^3$, there is no twisting surface scrolls of a minimal curve class. 
\end{remark}

\begin{lemma} \label{pretwisting} $X$ admits a very twisting maximal scroll if there exists an $1$-morphism $\zeta: \P^1\rightarrow \M_{0,1}(X, \theta)$ such that
	\begin{enumerate}
		\item the sheaf $\zeta^\ast T_{ev}$ is ample;
		\item the image $\zeta(\P^1)$ is in the smooth locus of the forgetful map $\Phi$ and the line bundle $\zeta^*T_\Phi$ is globally generated.
	\end{enumerate}
\end{lemma}
\proof Since $X$ is convex, every rational curve on $X$ is free. In particular, $(ev\circ \zeta)^*TX$ is globally generated. \qed

We may assume that $X$ is a projective homogeneous space under a connected semisimple linear algebraic $k$-group $G$. Let $T\subset G$ be a maximal torus. 

Let $\G_m\subset T$ be a one-dimensional torus corresponding to an interior point of a Weyl chamber. We recall basic properties of Bialynicki-Birula decompositions of $X$ under the torus action. See \cite{KimRahul}, \cite{BBdecomp}. The fixed points under the torus action  are isolated. For each $p\in X^{G_m}$, let $A_p$ be the set of points $x\in X$ such that $\lim_{t\rightarrow 0}t \cdot x=p$. By \cite[Prop. 1]{KimRahul}, $A_p$ is isomorphic to the affine space $\C^{l(p)}$, where $l(p)$ is the number of positive weights of the $\G_m$-representation at $T_p X$. 

Let $s, x_1,\dots,x_r\in X^{\G_m}$ be the fixed points corresponding to the unique maximal dimensional stratum $A_s$ and the set of all codimension one strata, $A_1,\dots,A_r$ respectively. Let $U$ be the union of $A_1,\dots,A_r$ and $A_s$, which is a dense open of $X$ with the complement at least codimension two. 

If we take the inverse torus action on $X$, there exists $1$-dimensional strata $A_1^\prime,\cdots, A_r^\prime$ corresponding to the fixed point $x_1,\dots,x_r$. Let $P_i$ be the closure of $A_i$, which is a smooth $\G_m$-invariant rational curve connecting $s$ and $x_i$. We call $P_i$'s the \emph{standard lines} on $G/P$ with respect to the $\G_m$-action. By \cite{KimRahul}, they generate the cone of effective curve classes of $G/P$.

\begin{lemma}\label{standardlineunique}
	The curve $P_i$ is the unique $\G_m$-invariant curve connecting $s$ and $x_i$. 
\end{lemma}

\proof By \cite[Prop. 1]{KimRahul}, there exists a $\G_m$-invariant open subset of $X$ containing $x_i$ which is $\G_m$-equivalent to a definite vector space representation $V_i$ such that the positive weight subspace of $V_i$ is of codimension one. Thus $P_i$ is the closure of the unique $G_m$-invariant curve in $V_i$ whose general point intersects $A_s$. \qed

\begin{definition}\label{transversaltype} Fix a $\G_m$-action on $X$ as above. A pointed maximal stable rational curve $f:(C,t_0)\rightarrow X$ is \emph{transversal}, if it satisfies the following properties:
	\begin{enumerate} 
		\item The image of $f(C)$ lies in $U$.
		\item The curve intersects $A_i$ transversally at $f(t_i)$.
		\item The marked point $f(t_0)$ is in $A_s$.
	\end{enumerate}\end{definition}
	
	A transversal maximal pointed rational curve $f$ gives an $(r+1)$-pointed rational curve $C^\prime=(C, t_0, t_1,\dots,t_r)$. 
	
	\begin{prop}\label{stablelimit}
		Given a transversal pointed maximal stable curve $f$ in $X$, the limit $\lim_{t\rightarrow 0} t\cdot f$ in $\M_{0,1}(X,\theta)$ is a $\G_m$-invariant pointed maximal stable rational curve $f_0: (F, p)\rightarrow X$ such that 
		\begin{enumerate}
			\item $F$ is obtained by gluing $\P^1_i$'s along the markings $t_i$'s of $C^\prime$, for $i=1,\cdots,r$,
			\item The marking $p$ is the point $t_0$ on $C^\prime$,
			\item the map $f_0$ maps $\P_i^1$'s to $P_i$ and contracts $C^\prime$ to $x_s$. 
		\end{enumerate}
	\end{prop} 
	
	\proof By Proposition \ref{simpleisunob}, $\M_{0,1}(X,\theta)$ is a smooth projective variety. Thus the limit under the torus action exists without the semistable reduction. The rest follows from \cite[Prop. 2]{KimRahul}.\qed
	
	There exists a natural map,
	$$\epsilon:\M_{0,1+r}\rightarrow \M_{0,1}(X,\theta)$$
	constructed as above. In fact, the morphsim $\epsilon$ is an isomorphism to its image by \cite{KimRahul}. The $\G_m$-action on $X$ induces the $\G_m$-action on $\M_{0,1}(X,\theta)$. By \cite{BBdecomp} and Proposition \ref{simpleisunob}, we consider the Bialynicki-Birula decomposition under the $\G_m$-action on $\M_{0,1}(X,\theta)$.
	
	\begin{cor}\label{bigstratum}
		Let $B$ be the image $\epsilon(\M_{0,1+r})$. The fixed locus $B$ is a smooth irreducible component of the $\G_m$-fixed point set in $\M_{0,1}(X,\theta)$ and the Bialynicki-Birula stratum corresponding to $B$ is of maximal dimensional.
	\end{cor}
	
	\proof The smoothness of $B$ is proved in \cite[Thm. 2.1]{BBdecomp}. A general maximal curve in $\M_{0,1}(X,\theta)$ is transversal by Kleiman-Bertini Theorem. By Proposition \ref{stablelimit}, it retracts to $\epsilon(\M_{0,1+r})$ under the $\G_m$-action. Thus there exists a dense open $\G_m$-invariant subset of $\M_{0,1}(X,\theta)$ retracting to the fixed point locus $B$, which by definition lies in the Bialynicki-Birula stratum of $B$.\qed
	
	\begin{lemma}\label{trivialtwisting}
		There exists an embedded rational curve in the fixed component $B$ such that the pullback of $T_\Phi$ and the normal bundle are positive. 
	\end{lemma}
	
	\proof With the discussion as above, the morphism $\epsilon:\M_{0,1+r}\rightarrow B$ is an isomorphism. Consider the forgetful map $F_0:\M_{0,1+r}\rightarrow \M_{0,r}$ by forgetting the first marked point. The fibers of $F_0$ give free curves in $\M_{0,1+r}$ such that the pullback of $T_\Phi$ is ample. We can choose a very free curve in $\M_{0,r}$ and lift it to a rational curve $D$ in $\M_{0,1+r}$. After attaching sufficiently many fibered curves of $F_0$ to $D$, a general smoothing of the comb yields the desired property. \qed

	Now we consider the inverse $\G_m$-action on $\M_{0,1}(X,\theta)$. By Corollary \ref{bigstratum}, There exists a fixed point component $B^\prime$ whose Bialynicki-Birula stratum is of maximal dimension. 
	
	Let $f: (C, p)\rightarrow X$ be a general maximal rational curve in $X$. We may assume that $[f]$ lies in both Bialynicki-Birula strata corresponding to $B$ and $B^\prime$. Let $\zeta:\P^1\rightarrow \M_{0,1}(X,\theta)$ be a $\G_m$-orbit curve of $[f]$. The image $\zeta(0)$, resp., $\zeta(\infty)$ corresponds to a $\G_m$-invariant curve $[f_0]$ in $B$, resp., $[f_\infty]$ in $B^\prime$. By \cite[Thm. 4.3]{BBdecomp}, we have the following $\G_m$-equivariant decomposition of the tangent spaces,
	$$T_{[f_0]}\M_{0,1}(X,\theta)=T_{[f_0]}B \oplus T_{[f_0]}\M_{0,1}(X,\theta)^+,$$
	$$T_{[f_\infty]}\M_{0,1}(X,\theta)=T_{[f_\infty]}B^\prime \oplus T_{[f_\infty]}\M_{0,1}(X,\theta)^-.$$
	Here the $\G_m$-actions on $T_{[f_0]}B$ and $T_{[f_\infty]}B$ are both trivial and $T_{[f_0]}\M_{0,1}(X,\theta)^+$ ($T_{[f_\infty]}\M_{0,1}(X,\theta)^-$) corresponds to the positive (negative) weight $\G_m$-invariant subspace. Since the evaluation map $ev:\M_{0,1}(X,\theta)\rightarrow X$ is $\G_m$-equivariant and smooth, we have the sub-decompositions of $T_{ev}$: 
	$$T_{ev,[f_0]}=T_{[f_0]}B \oplus T_{ev,[f_0]}^+,$$
	$$T_{ev,[f_\infty]}=T_{[f_\infty]}B^\prime \oplus T_{ev,[f_\infty]}^-.$$
	The decomposition of weight spaces at $T_{ev,[f_0]}$ uniquely determines a decomposition of the $\G_m$-equivariant vector bundle $\zeta^*T_{ev}$, i.e., $$\zeta^*T_{ev}=E^0\oplus E^+,$$ 
	where $E^0|_{[f_0]}=T_{[f_0]}B$ and $E^+|_{[f_0]}=T_{ev,[f_0]}^+$. 
	
	\begin{prop} \label{1-twisting}A general $\G_m$-orbit curve $\zeta:\P^1\rightarrow \M_{0,1}(X,\theta)$ satisfies the following:\begin{enumerate}
			\item The sheaf $E^0$ is a semi-positive vector bundle over $\P^1$. 
			\item The sheaf $E^+$ is a positive vector bundle over $\P^1$.
			\item The image $\zeta(\P^1)$ is in the smooth locus of $\Phi$ when $r\neq 2$. The line bundle $\zeta^*T_{\Phi}$ is positive when $r=1$, and is trivial when $r\ge 3$.
		\end{enumerate}
		
	\end{prop}
	
	\proof By the definition of $E^0$ and $E^+$ as above, the weights of $E^0$, resp., $E^+$ at $0$ are trivial, resp., positive. The weights of $E^0$ and $E^+$ at $\infty$ are both non-positive. Since the degree of any $\G_m$-equivariant line bundle equals the difference of the weight at $0$ and the weight at $\infty$, we get (1) and (2).
	
	For (3), note that $\zeta^*T_\Phi$ is a $\G_m$-equivariant vector bundle on $\P^1$. When $r=1$, the curve $[f_0]$ is a pointed line $L$ in $X$ by Proposition \ref{stablelimit}. Thus $T_{\Phi,[f_0]}$ is isomorphic to $T_p L$ as a vector space. The weight is positive because the marked point is a retracting fixed point. Similarly, the weight at $T_{\Phi,[f_\infty]}$ is negative. Hence $\zeta^*T_{\Phi}$ is a positive line bundle.
	
	When $r\ge 3$, the marked point on $[f_0]$, resp. $[f_\infty]$ lies in the contracted component and as well as in the smooth locus of $\Phi$. Thus the weight at $0$ and $\infty$ are both trivial under the torus action, i.e., $\zeta^*T_\Phi$ is a trivial vector bundle.\qed
	
	\begin{prop}\label{pic12}
		When the Picard number of the homogeneous space $X$ is either one or two, there exists a very twisting maximal scroll on $X$.
	\end{prop}
	
	\proof With the same notations as above, in either case, the fixed locus $B$ which corresponds to the maximal Bialynicki-Birula cell is a point. Hence, as in Proposition \ref{1-twisting}, for a general $\G_m$-orbit curve $\zeta$, there is no $E^0$-summand in $T_{ev}$. Thus the weights of the $\G_m$-vector bundle $\zeta^*T_{ev}$ at $0$, resp., at $\infty$, are all positive, resp., negative. Therefore, $\zeta^*T_{ev}$ decomposes into a direct sum of line bundles with degrees $\ge 2$. 
	
	When the Picard number is one, by Lemma \ref{pretwisting} and the third part in Proposition \ref{1-twisting}, we win.
	
	When the Picard number of $X$ is two, we have trouble analyzing $T_\Phi$ because the two $\G_m$-fixed points $\zeta(0)$ and $\zeta(\infty)$ lie in the singular locus of $\Phi$. However, the singular locus of $\Phi$ in $\M_{0,1}(X,\theta)$ is of codimension two. Note that the orbit curve $\zeta$ is free in $\M_{0,1}(X,\theta)$. Hence, a general deformation $\xi: \P^1\rightarrow \M_{0,1}(X,\theta)$ of $\zeta$ avoids the singular locus of $\Phi$ and intersects the boundary divisors of $\M_{0,1}(X,\theta)$ transversally. The pullback of the universal family over $\M_{0,1}(X,\theta)$ over $\xi$ gives a smooth surface $S$ over $\P^1$ with a section $D$. The sheaf $\xi^*T_{ev}$ is positive by upper semicontinuity. The degree of the line bundle $\xi^*T_\Phi$ is the self-intersection number $(D.D)$ on $S$, which is constant in the deformed family. Thus it suffices to check for $\zeta$. The marked point in universal family over $\zeta$ gives a section in the smooth locus with self-intersection zero. See \cite[Prop. 2]{KimRahul}. In particular, $\xi^* T_\phi$ is trivial. By Lemma \ref{pretwisting}, a general deformation of $\zeta$ gives a very twisting maximal scroll on $X$.\qed

	To construct a very twisting surface maximal scroll on projective homogeneous space of higher Picard numbers, the main idea is to glue a bunch of ``nearly" very twisting scrolls as above properly whose general smoothing is very twisting.
	
	\begin{construction}\label{cons} Let $X$ be projective homogeneous spaces with the Picard number greater than two. The $\G_m$-fixed component $B$ in (\ref{bigstratum}) has positive dimension. By Lemma \ref{trivialtwisting}, there exists a rational curve $D$ in $B$ such that both $\N_{D|B}$ and $T_\Phi|_D$ are postive vector bundles. Since $D$ is very free, we may choose distinct points $p_1,\cdots,p_k$ on $D$, where $p_i$ is the limit point of a $\G_m$-orbit curve $C_i$ as in Proposition \ref{1-twisting}. Let $C$ be the disjoint union $\coprod_{i=1}^k C_i$. Consider the comb $D^*=D+\sum_{i=1}^k C_i=D+C$ obtained by attaching each $\G_m$-orbit curve $C_i$ on $D$ at $p_i$. \end{construction}
	

	
	\begin{lemma}\label{smoothing}
		After attaching sufficiently many general $C_i$'s on $D$, the comb $D^*$ can be smoothed. 
	\end{lemma}
	
	\proof By \cite[Lem. 2.6]{GHS}, the normal sheaf $\N_{D^*}$ restricted on $D$ is the sheaf of rational sections of $\N_D$ having at most a simple pole at each $p_i$ in the normal direction determined by $T_{p_i}C_i$. By the short exact sequence,
	
	$$\begin{CD}0@>>>\N_{D|B}@>>>\N_D@>>>\N_B|_D@>>>0\end{CD}$$
	the normal directions in $\N_D$ determined by $T_{p_i}C_i$'s give nonzero general directions in $\N_B|_D$. Thus the quotient bundle $\mathcal{M}=\N_{D^*}|_D/\N_{D|B}$ is the sheaf of rational sections of $\N_B|_D$ having at most a simple pole at each $p_i$ in the normal direction determined by $T_{p_i}C_i$. By \cite[Lem. 2.5]{GHS}, after attaching sufficiently many general $C_i$'s, $\mathcal{M}$ is globally generated. Together with the positivity of $\N_{D|B}$, the sheaf $\N_{D^*}|_D$ is globally generated. Since all $C_i$'s are free, by diagram chasing, the normal sheaf $\N_{D^*}$ is globally generated. In particular, the comb $D^*$ is unobstructed and the nodes can be smoothed. \qed
	
	Choose a smoothing of $D^*$ over a smooth pointed curve $(T,0)$ as the following,
	
	
	$$\xymatrix{
		D^* \ar[d] \ar@{^{(}->}[r]  &
		S   \ar[d]^p \ar@{^{(}->}[r] & \M_{0,1}(X,\theta)\\
		0 \ar[r]& (T,0) }$$
	where $S$ is a smooth surface. Let $\E$ be the pullback bundle of $T_{ev}$ to $S$. Let $E^0_i$, resp., $E^+_i$ be the trivial, resp., positive subbundle of $T_{ev}$ restricted to each $C_i$. Let $\T$ be the vector bundle $\coprod E_i^+$ over $C$. Since $\T$ is a direct summand of $\E|_C$, we have the following natural surjection. $$\E^\vee\rightarrow \E^\vee|_{C}\rightarrow \T^\vee$$
	Let $\K^\vee$ be the elementary transform of $\E^\vee$ along $\T^\vee$. 
	
	\begin{equation}\label{dual}\begin{CD}0@>>>\K^\vee @>>>\E^\vee@>>>\T^\vee @>>>0\end{CD}\end{equation}
	Dualizing the above short exact sequence, we get
	\begin{equation}\label{ndual}\begin{CD}0@>>>\E @>>>\K @>>>\T\otimes_{\O_C}\O_C(C) @>>>0\end{CD}.\end{equation}

	\begin{lemma}\label{posquill}
		For any $i=1,\cdots,k$, $h^1(C_i, \K|_{C_i}(-p_i))=0$.
	\end{lemma}
	
	\proof Restricting the short exact sequence (\ref{dual}) to $C_i$ and applying the functor $\Hom_{\O_{C_i}}(\  \_\_\  ,\O_{C_i})$, we get the following exact sequence
	$$\begin{CD}0@>>>E_i^+ @>>>\E|_{C_i} @>>>\K|_{C_i}@>>>E_i^+\otimes_{\O_{C_i}}\O_{C_i}({C_i})@>>>0\end{CD}.$$
	The quotient bundle $\E|_{C_i}/E_i^+$ is $E_i^0$ and the last term of the exact sequence is isomorphic to $E_i^+(-p_i)$. In particular, we have 
	$$\begin{CD}0@>>>E_i^0(-p_i)  @>>>\K|_{C_i}(-p_i)@>>>E_i^+(-2p_i)@>>>0\end{CD}.$$
	Note that over $C_i$, $E_i^0$ is trivial and $E_i^+$ is positive. We win.\qed 
	
	Let $s_1$ and $s_2$ be two sections of $p$ both of which specialize to two distinct point $q_1, q_2$ on $D^*\backslash C$. 
	
	\begin{lemma}\label{posbody} We have $h^1(D,\K|_{D}(-p_1-p_2))=0$, after attaching sufficiently many $C_i$'s on $D$.
	\end{lemma}
	
	\proof Restricting the short exact sequence (\ref{dual}) to $D$, we get 
	$$\begin{CD}\K^\vee|_{D} @>>>\E^\vee|_{D}@>>>T^\vee|_D@>>>0\end{CD}.$$
	The above sequence is actually exact. Indeed, by restricting (\ref{ndual}) to $D$ and taking the dual over $D$, since $\T\otimes_{\O_C}\O_C(C)|_D$ is torsion, we have the injection from $\K^\vee|_{D}$ to $\E^\vee|_{D}$. 
	
	In other words, the vector bundle $\K^\vee|_D$ is the elementary transform up of $\E|_D$ along $p_i$'s with the specific directions in $E_i^+$'s. Since the sub-bundle $TB|_D$ of $\E|_D$ restricting to each $p_i$ is orthgonal to $\T|_{p_i}=E_i^+$, it is also a sub-bundle of $\K|_D$.
	
	Since $TB|_D$ is ample, to prove the Lemma, it suffices to show that the quotient bundle $(\K|_D)/(TB|_D)$ is positive on $D$ after attaching sufficiently many $C_i$'s. Consider the following diagram.
	
	$$\begin{CD}
	0@>>>{(\frac{\K|_{D}}{TB|_D})}^\vee @>>>(\frac{\E|_{D}}{TB|_D})^\vee @>t>>T^\vee|_D@>>>0\\
	@. @VVV @VVV @| @.\\
	0@>>>\K^\vee|_{D} @>>>\E^\vee|_{D}@>>>T^\vee|_D@>>>0\\
	@. @VVV @VVV @VVV @.\\
	0@>>>(TB|_D)^\vee @= (TB|_D)^\vee @>>> 0 @.
	\end{CD}.$$
	We get that the vector bundle $(\K|_D)/(TB|_D)$ is the elementary transform up of $(\E|_D)/(TB|_D)$ along $p_i$'s with the direction $E_i^+$'s. Note that the torsion quotient $t$ is just the restriction of $(\E|_D)/(TB|_D)$ at $p_i$'s. Thus $(\K|_D)/(TB|_D)$ is isomorphic to $(\E|_D)/(TB|_D)\otimes_{\O_D}\O_D(\sum p_i)$, which is positive when the attachment points on $D$ are sufficiently many.\qed

	\begin{theorem}\label{twisting}
		Let $X$ be a projective homogeneous space over an algebraically closed field $k$ of characteristic zero. Let $\theta$ be the maximal curve class on $X$. There exists a very twisting maximal scroll $\zeta: \P^1\rightarrow \M_{0,1}(X, \beta)$.
	\end{theorem}


	
	\proof By Proposition \ref{pic12}, it suffices to prove the case when $X$ has Picard number greater than two. Now we may construct the comb $D^*$ as in (\ref{cons}) by attaching sufficiently many general $C_i$'s. By Lemma \ref{smoothing}, the comb can be smoothed. By Lemmas \ref{posquill} and \ref{posbody}, $h^1(D^*,T_{ev}|_{D^*}(-s_1-s_2))$ is zero. Thus by upper semi-continuity, $T_{ev}$ restricting to a general smoothing of $D^*$ is ample. 
	
	Similarly Condition (3) of Proposition \ref{1-twisting} and Lemma \ref{trivialtwisting}, the vector bundle $T_\Phi|_{D^*}$ is positive. Therefore $T_\Phi$ restricting to a general smoothing of the comb $D^*$ is also positive by upper semi-continuity. The theorem is proved by Lemma \ref{pretwisting}.\qed
	
	\section{Rational Simple Connectedness of Homogeneous Spaces}
	
	\begin{prop}\label{H1done}
		Let $X$ be a projective homogeneous space defined over an algebraically closed field of characteristic zero. Then for any simple curve class $\beta$, the evaluation morphism
		$$ev:\M_{0,1}(X,\beta)\rightarrow X$$
		is smooth surjective with integral rationally connected geometric fibers.
	\end{prop}
	
	\proof The evaluation map $ev$ is smooth because of the generic smoothness and the homogeneity of the target $X$. Since $X$ is simply connected, the finite part of the Stein factorization of $ev$ is \'etale over $X$, thus isomorphic to $X$. Therefore every geometric fiber is connected and smooth, thus integral.
	
	By Proposition \ref{simpleisunob}, the moduli space $\M_{0,1}(X,\beta)$ is a nonempty smooth projective rational variety. By \cite[Lem. 15.6]{dJHS}, the geometric fibers of the evaluation morphism are rationally connected.\qed
	
	Let $k$ be an algebraically closed field of characteristic zero. Let $G$ be a connected reductive linear algebraic group over $k$. Let $T\subset G$ be a maximal torus of rank $t$ and let $B$ be a Borel subgroup of $G$ containing $T$. The choice of $(G,B,T)$ gives a root system. Let $\Delta=\{\alpha_1,\cdots,\alpha_t\}$ be a basis of the root system. Let $W$ be the Weyl group of the root system generated by simple reflections $\{s_i=s_{\alpha_i}|\alpha_i\in\Delta\}$.
	
	Let $n_w\in N_G(T)$ be a representative of $w\in W$. The map $w\mapsto n_w B$ induces a one-to-one correspondence between the Weyl group and the set of $T$-fixed points in $G/B$. We simply write $w$ for the corresponding fixed point.
	
	Let $U$ be the unipotent radical of $B$. By Bruhat decomposition \cite[14.12]{Borel}, $G/B$ is a disjoint union of $U$-orbits $Uw$ and each orbit is isomorphic to the vector space $k^{l(w)}$, where $l$ is the length function on the Weyl group. Let $w_0$ be the longest element of $W$. It corresponds to the maximal dimensional Bruhat cell. Let $w_1,\cdots,w_t$ be the fixed points of $G/B$ which correspond to the codimension one Bruhat cells. 
	
	Let $\G_m\subset T$ correspond to the interior of the positive Weyl chamber. By \cite[3.4.7]{Carrell1}, the Bialynicki-Birula decomposition of $G/B$ coincides with the Bruhat decomposition. Thus each standard line in $G/B$ is the unique $\G_m$-invariant line connecting $w_0$ and $w_i$. 
	
	\begin{lemma}\label{algequiv}
		Every maximal curve in $G/B$ is algebraically equivalent to the union of all standard lines.
	\end{lemma}
	
	\proof This is a corollary of Proposition \ref{simpleisunob} and Proposition \ref{stablelimit}.\qed
	
	Let $I$ be a subset of $\Delta$. Let $W_I$ be the subgroup of the Weyl group generated by simple reflections of $I$. The \emph{standard parabolic subgroup} is of the form $BW_IB$. Every parabolic subgroup of $G$ is conjugate to the standard parabolic subgroup $P_I$ containing $B$. Thus every projective homogeneous space under $G$ is of the form $G/{P_I}$.
	
	Let $\pi_I:G/B\rightarrow G/{P_I}$ be the natural projection. The induced $\G_m$-action on $G/P_I$ induces a one-to-one correspondence between the $\G_m$-fixed points and the left coset space $W/W_I$. For each coset $wW_I$, there exists a unique representative $w'$ with the minimal length and  and $l(w'w'')=l(w')+l(w'')$ for any $w''\in W_I$, cf. \cite[1.10]{Humphreys1}. By \cite[3.4.8]{Carrell1}, each Bialynicki-Birula cell of $wW_I$ is isomorphic to $k^{l(w')}$. It is easy to see that $w_0=w^0w_{I_0}$, where $w_{I_0}$ is the longest element in $W_I$ and $l(w^0)$ is the dimension of $G/P$.
	
	\begin{lemma} \label{standardlift} For each standard line in $G/P_I$, there exists a unique lifting to a standard line in $G/B$.
	\end{lemma}
	
	\proof First we show that every fixed point in $G/P$ corresponding to a codimension one cell uniquely lifts to a fixed point in $G/B$ satisfying the same property. For each coset $wW_I$ with the representative $w'$ discussed above, $w'w_{I_0}$ is the unique element in $wW_I$ with maximal length. If a coset $wW$ corresponds to a codimension one cell in $G/P$, i.e., $l(w')=l(w^0)-1$, we have 
	$$l(w'w_{I_0})=l(w')+l(w_{I_0})=l(w^0)+l(w_{I_0})-1=l(w_0)-1.$$
	Thus the fixed point $w'w_{I_0}$ in $G/B$ corresponds to a unique codimension one cell. 
	
	The standard line $L$ connecting $w_0$ and $w'w_{I_0}$ in $G/B$ projects to a $\G_m$-invariant curve connecting $w_0W$ and $w'W$ in $G/P$. By Lemma \ref{standardlineunique}, the image $\pi_I(L)$ is a standard line in $G/P$. Since the projection morphism between the big cell of $G/B$ and the big cell of $G/P$ is a $\G_m$-equivariant linear morphism between vector spaces, the degree of $\pi_I|_L$ is one. Thus $L$ maps isomorphically onto its image, which is a standard line. We get the lifting.\qed 
	

	
	
	\begin{lemma}\label{maximalsimple}
		Every maximal curve in $P_I/B$ gives a simple curve of $G/B$.
	\end{lemma}
	
	\proof With the $\G_m$-action on $G/B$ as above, by Lemma \ref{algequiv}, it suffices to show that standard lines in $P_I/B$ correspond to standard lines in $G/B$ and the correspondence is injective. Any standard line in $P_I/B$ is the unique $\G_m$-invariant line connecting $w_{I_0}$ and $w_{I_0}s_i$, where $t_i\in I$ by Lemma \ref{standardlineunique}. After the left translation by $w^0$, we get a $\G_m$-invariant line connecting $w_0$ and $w_0s_i$, which is standard in $G/B$ by Lemma \ref{standardlineunique} again. Since such correspondence is induced by a left translation, clearly it is injective. \qed
	
	\begin{prop}[\cite{dJHS}, Def. 7.1]\label{chainsmooth}
		The moduli space $\Chn_2(X,m\theta)$ of two-pointed chains of $m$ stable maximal curves in $X$ is represented by a nonempty smooth projective variety.\qed
	\end{prop}

	\begin{prop}\label{H2done}
		Let $X$ be a projective homogeneous space defined over an algebraically closed field of characteristic zero. Then there exists $m$ such that the geometric generic fiber of the evaluation morphism
		$$ev:\Chn_{2}(X,m\theta)\rightarrow X\times X$$
		is smooth integral rationally connected.
	\end{prop}
	
	\proof By Corollary \ref{chainsmooth}, the moduli space of two-pointed chains of $m$ maximal curves is a smooth projective variety. By induction on $m$ and Proposition \ref{simpleisunob}, it is rationally connected. By the proof of \cite[Lem. 15.8]{dJHS}, it suffices to show that the evaluation $$ev:\Chn_{2}(X,m_0\theta)\rightarrow X\times X$$
	is surjective for some $m_0$. Assume that $X=G/P$, where $G$ is a reductive group. We prove this by induction on the rank of $G$. By Lemma \ref{algequiv} and Lemma \ref{standardlift}, 
	it suffices to show the case when $X=G/B$. When the rank of $G$ is one, the surjectivity of $ev$ is trivial because $G/B$ is isomorphic to $\P^1$.
	
	When the rank of $G$ is bigger than one, let $\Delta$ be the set of simple roots of $G$. Let $P_i$ be the standard parabolic subgroup corresponding to a simple root $\alpha_i\in\Delta$. Let $P^i$ be the standard maximal parabolic subgroup corresponding to $\Delta-{\alpha_i}$. Let $s_i$ be the simple reflection of $\alpha_i$. Consider the following diagram,
	$$\begin{CD}
	G/B@>u>>G/P^i\\
	@VvVV\\
	G/P_i
	\end{CD}$$
	where $G/P^i$ is a projective homogeneous space of Picard number one and the morphism $v$ is a $\P^1$-bundle over $G/P_i$. By the proof of Lemma \ref{maximalsimple}, the fiber of $v$ is algebraically equivalent to the standard line $L_i$ through $w_0$ and $w_0s_i$ in $G/B$. Since $s_i$ is not in $W_{\Delta-\{\alpha_i\}}$, the images $u(w_0)$ and $u(w_0s_i)$ are disjoint in $G/P^i$. By Lemma \ref{standardlift}, $L_i$ maps to the unique standard line in $G/P^i$. Thus all the fibers of $v$ map to lines in $G/P^i$. We call the image lines in $G/P^i$ \emph{good} lines. In fact, the above diagram gives a connected proper flat prerelation on $G/P^i$. By \cite[IV.4.14]{Kollar} and by homogeneity, every pair of points in $G/P^i$ can be connected by a chain of good lines of length $m$. 
	
	Now given a pair of points $p$ and $q$ in $G/B$,  there exists a chain of $m$ good lines in $G/P^i$ connecting $u(p)$ and $u(q)$. We can lift the good lines to $m$ two pointed lines $(l_1,p_1,q_1),\cdots (l_m,p_m,q_m)$ in $G/B$ such that $u(p_1)=u(p)$, $u(q_m)=u(q)$, and $u(q_{i})=u(p_{i+1})$ for $i=1,\cdots, m-1$. 
	
	The fiber of $u$ is a projective homogeneous space under an algebraic group of smaller rank, i.e., a Levi subgroup of $P_i$. By induction, we can choose chains of maximal curves in the fiber of $u$, connecting $p$ and $p_1$, $q_1$ and $p_2$, etc. By Lemma \ref{maximalsimple}, we get a chain of simple curves in $X$ connecting $p$ and $q$. By adding lines to make each irreducible component of the chain maximal, we get a maximal chain connecting $p$ and $q$ in $G/B$.\qed


\section{On Discriminant Avoidance}
Let $k$ be an algebraically closed field of arbitrary characteristic. Let $S$ be a $k$-variety of dimension $d$. Let $K$ be the function field of $S$. Let $X$ be a smooth projective Fano $k$-variety and $U$ be its universal torsor over $X$. Let $r$ be the Picard number of $X$. Since $k$ is algebraically closed, $U$ is a $(\G_m)^r$-torsor over $X$ and $U$ exists unique up to isomorphism. We consider the following question. 

\begin{question}
	Given $p:\X\rightarrow S$ an isotrivial family of $X$ over $S$ with the vanishing of the elementary obstruction on the generic fiber, is there a rational section? 
\end{question}

By Proposition \ref{torsor-ob}, the vanishing of the elementary obstruction is equivalent to the existence of the universal torsor of $\X_K$. After shrinking the base $S$ to an open subset, the above question is equivalent to the following.

\begin{question}\label{keyQ}
	Given $(p:\X\rightarrow S, \U)$ an isotrivial family of $(X,U)$ over $k$, is there a rational section? 
\end{question}

Let $G$ be the automorphism group of the pair $(X,U)$ over $k$. The group scheme $G$ has $T$-valued points which are the pairs $(\phi,\alpha)$, where $\phi: X_T\rightarrow X_T$ is an automorphism of schemes over $T$ and $\alpha:\phi^*U\rightarrow U$ is an isomorphism of $(\G_m)^r$-torsors. 

The question \ref{keyQ} gives $(p:\X\rightarrow S, \U)$, which is an isotrivial family of the pair $(X,U)$ over $S$. It is natural to associate the pair with a $G$-torsor over $S$. Consider the functor that the $T$-valued points over $S$ are the set of pairs $(\phi,\alpha)$, where $\phi: \X_T\rightarrow \X_T$ is an automorphsim of schemes over $T$ and $\alpha:\phi^*\U\rightarrow \U$ is an isomorphism of $Hom_{T}(R^1p_{T*}\G_m,\G_{m,T})$-torsors. 

\begin{lemma}
	If $S$ is reduced, the functor is representable by a scheme $\T$ over $S$ and $\T$ is a $G$-torsor over $S$ by post-composing.
\end{lemma}

\proof Since every $G$-torsor over $S$ is affine, it suffices to prove the representability of the functor fppf locally by the descent of affine group schemes. First we will show that the pair $(p:\X\rightarrow S, \U)$ is fppf locally isomorphic to the constant family. 

By taking an \'etale neighborhood $V$, we may assume that the pullback of the torsor $\U$ is a $\G_m^r$-torsor over $\X_V$. Thus the relative character lattice is isomorphic to $\Z^r\times V$. We can choose a basis $L_1,\dots,L_r$ of the relative character lattice such that each $L_i$ corresponds to a very ample line bundle ($\G_m$-torsor) over $\X|_V$. Now by the Hilbert scheme trick used in the proof of Lemma 2.2.1 in \cite{dJS2}, after a flat base change, the pairs $(\X|_V, L_i)$ are constant families. So is the pair $(\X|_V,\U|_V)$.

This implies that the functor restricted on $V$ is just $\Isom_V((X_V, U_V),(X_V, U_V))$ and $U_V$ is a $(\G_m)^r$-torsor over $X_V$. Since $X$ is Fano, we know that $\Aut(X)$ is represented by a linear algebraic group. Thus $\Isom_V((X_V, U_V),(X_V, U_V))$ is represented by the scheme $G\times V$. This proves the lemma.\qed

\begin{lemma}
	Given a $G$-torsor $\T$ over $S$, we can associate a pair $(p:\X\rightarrow S, \U)$ where $\U$ is a relative universal torsor over $\X$.
\end{lemma}

\proof The morphism $\T\rightarrow S$ is fppf. It suffices to descent the constant family $(X,U)\times \T$ to $S$. First we will descent the isotrivial family of $X$. Since such family has a natural polarization, the anti-canonical polarization, it is easy to check that the polarized family descents to $S$. Similarly, we can descent the relative Picard scheme and the torsor under the relative Picard scheme to $S$ by \cite[Ch. 5, Sec. 6]{BLR}. The new torsor being universal follows from the universality of the constant family, cf., \cite[Prop. 2.2.4]{SK}. \qed

\begin{theorem} 
	If $G=\Aut(X,U)$ is geometrically reductive, then Question \ref{keyQ} can be reduced to the projective base case. 
\end{theorem}

\begin{remark}
	This is called \emph{discriminant avoidance}, which is studied by de Jong and Starr \cite{dJS2} for isotrivial families of Picard number $1$. For varieties of higher Picard numbers, it is natural to replace ample generating line bundles in their setting by universal torsors. The latter gives a cohomological obstruction to the existence of rational points.
\end{remark}

\proof
By the above two lemmas, we get a one-to-one correspondence between isotrivial families $(p:\X\rightarrow S, \U)$ and $G$-torsors over $S$ when $S$ is reduced. The remaining part is exactly the same as the proof of Theorem 2.1.3 in \cite{dJS2}.\qed

The following Lemma gives a description of $G=\Aut(X,U)$.
\begin{lemma}
	If $X$ is Fano, then $G=\Aut(X,U)$ is an extension of $\G_m^r$ and $\Aut(X)$, where $\Aut(X)$ is a linear algebraic group. In particular, if $\Aut(X)$ is geometrically reductive, $G$ is geometrically reductive.
\end{lemma}

\proof Since $X$ is Fano, we can choose a large multiple of the anticanonical bundle to embed $X$ into a projective space. Thus $\Aut(X)$ is a linear subgroup of $PGL(N)$. There is a left exact sequence of linear algebraic groups, where $\Aut_X(X,U)$ is the kernel of the forgetful map.
$$\begin{CD}
1@>>>\Aut_X(X,U)@>>>\Aut(X,U)@>F>>\Aut(X)
\end{CD}$$
By \cite[Lem. 4.1]{Brion}, $\Aut_X(X,U)$ is isomorphic to the group $\Hom(X,\G_m^r)$. Since $X$ is projective,  $\Hom(X,\G_m^r)\cong\G_m^r$. 

It suffices to show that the forgetful map $F$ is surjective. For any automorphism $\phi$ of $X$, the pullback $\phi^*U$ is again a universal torsor. The universal torsor is unique up to isomorphism over $X$ when $k$ is algebraically closed. We can choose any isomorphism between $\phi^*U$ and $U$. \qed

\begin{cor}\label{DA}
	The discriminant avoidance holds for isotrivial families of Fano varieties if the automorphism group of the fiber is geometrically reductive.\qed
\end{cor}

\section{Proof of the Main Theorem}

\begin{lemma}\label{picardreduction}
	Let $X$ be a projective homogeneous space defined over a field $K$. Assume that the elementary obstruction vanishes and the Picard number of $X$ is greater than one. Then there exists a smooth morphism,
	$$\begin{CD}X@>u>>Y@>>>\Spec K\end{CD}$$
	such that $Y$ is a projective homogeneous space of Picard number one with the vanishing elementary obstruction. Furthermore, if $Y$ admits a rational point $p$, then the fiber $u^{-1}(p)$ is a smooth projective homogeneous space with the vanishing elementary obstruction.
\end{lemma}

\proof Let $\Gamma$ be the Galois group of the field $K$. When the elementary obstruction of $X$ vanishes, by \cite[Prop. 2.25]{CTS}, $\Pic(X)$ is isomorphic to $\Pic(\overline{X})^\Gamma$. Thus by assumption the rank of $\Pic(\overline{X})^\Gamma$ is greater than one. By Lemma \ref{lemperm}, $\Pic(\overline{X})$ is a permutation $\Gamma$-module with a canonical $\Gamma$-invariant basis $\L_1,\cdots, \L_r$. We can choose a $\Gamma$-orbit in the basis, denoted by $\L_1,\cdots,\L_b$. Since $\L=\L_1+\cdots+\L_b$ is $\Gamma$-invariant, the line bundle $\L$ is globally generated and defined over $K$. The linear system $|\L|$ gives the morphism $u:X\rightarrow Y$. It is clear from the construction that $u$ is smooth and $Y$ is a projective homogeneous space and of Picard number one. The vanishing of the elementary obstruction of $Y$ follows from \cite[Lem. 3.1.2]{Wittenberg0}.

Let $\overline{X}$ be the base change of $X$ to the algebraic closure. A universal torsor on $\overline{X}$ is isomorphic to a $\G_m^r$-torsor $\L_1 \times\cdots\L_r$ which is unique up to isomorphism. The vanishing of the elementary obstruction is equivalent to that the universal torsor on $\overline{X}$ descents to $X$, cf., \cite[Prop. 2.2.4]{SK}. Let $\T$ be the universal torsor on $X$ and $\T_p$ be the restriction of $\T$ on $Z=u^{-1}(p)$. By functorality of the restriction, $\T_p\times_{K} \overline{K}$ is the same as $\T\times_{K} \overline{K}|_{\overline{Z}}$. The latter term is just $\L_1 \times\cdots\L_r|_{\overline{Z}}$. It is easy to see that the restriction gives a product of a trivial $\G_m^b$-torsor and the universal torsor on $\overline{Z}$. Therefore the elementary obstruction of $Z$ vanishes.\qed

\begin{lemma}\label{auto}
	Let $X$ be a projective homogeneous space $G/P$ over an algebraically closed field of characteristic zero. Then the connected component of the automorphism group $\Aut(X)$ is reductive.
\end{lemma}

\proof Since $X$ is Fano, the automorphism group is a linear algebraic group. Let $R$ be the solvable radical of the connected component of $\Aut(X)$. The solvable group $R$ naturally acts on $X$. By the Borel fixed point theorem \cite[III.10.4]{Borel}, there exists a fixed point $x$ of $R$. Let $L_g$ be the automorphism of the left translation on $X$ by an element of $g\in G$, which clearly lies in the connected component of $\Aut(X)$. For any closed point $y$ in $X$, there exists $g\in G$ such that $L_g(y)=x$. For every element $\varphi$ in $R$, since $R$ is normal, $L_{g}\circ\varphi\circ L_{g^{-1}}$ lies in $R$. Thus we have
$$L_g(\varphi(y))=(L_{g}\circ\varphi\circ L_{g^{-1}})( L_g(y))=(L_{g}\circ\varphi\circ L_{g^{-1}})(x)=x=L_g(y).$$  
Thus $\varphi$ fixes $y$, i.e. $\varphi$ fixes every point in $X$. This implies that the solvable radical $R$ is trivial.\qed

\proof[Proof of Theorem \ref{intromain}] By Proposition \ref{abs}, we only need to prove the ``if" case. By \cite[Lem. 16.3]{dJHS}, it suffices to prove the theorem in characteristic zero. By Lemma \ref{picardreduction} and induction on the Picard number, it suffices to prove the case when the Picard number of $X$ is one. Let $\pi:\X\rightarrow U$ be an integral model of $X$, where $U$ is a dense open subset of $S$. After shrinking $U$, we may assume that $\pi$ is smooth and the relative universal torsor exists. By the method of discriminant avoidance, cf., Lemma \ref{auto} and Corollary \ref{DA}, we may assume that $U=S$ is projective.

After blowing up the base points of a Lefschetz pencil of $S$,
we have the right column of the following diagram. When taking the base change to the generic point of $\P^1$, we have the left column of the following Cartesian diagram.
$$\begin{CD}
X@>>>\X\\
@V{\pi}VV @VVV\\
C@>>>S\\
@VVV @VVV\\
k(\P^1)@>>> \P^1
\end{CD}$$

Let $K$ be the field $k(\P^1)$. Now we are in Situation \ref{hyp}. By Proposition \ref{H1done} and \ref{H2done}, Hypotheses \ref{H1} and \ref{H2} hold. By Theorem \ref{twisting}, Hypothesis \ref{H3} holds. By Theorem \ref{major}, there exists an Abel sequence $(Z_e)_{e\ge e_0}$ for $X/C/K$.

Therefore the Abel map $\alpha:Z_e\rightarrow \Pic_{D/K}^e$ is surjective with integral rationally connected geometric generic fiber for $e\gg 0$. Since the exceptional curves on $S$ give the constant sections of $S\rightarrow \P^1$, there exist rational points on $\Pic_{C/K}^e$ for every integer $e>0$. By pullback to $D$, there exist rational points on $\Pic_{D/K}^{re}$ for every $e>0$, where $r$ is the geometric Picard number of $X$. When $e\gg 0$ and divisible by $r$, the fiber of the Abel map over a rational point of $\Pic_{D/K}$ is integral rationally connected defined over $K$. By \cite{GHS}, there exists a $K$-rational point on the coarse moduli space of $Z_e$. By \cite[Lem. 13.3]{dJHS}, we get a rational point.\qed

\begin{lemma}[Starr]\label{jasontrick}
	Let $K$ be a field. Let $G$ be a quasisplit adjoint semisimple group defined over $K$. If a $G$-torsor admits a reduction to a Borel subgroup, then it is trivial. 
\end{lemma}

\proof  Let $\Won(G)$ be the wonderful compactification of $G$. For any $G$-torsor $\T$, we can twist $\Won(G)$ by $\T$ using the right $\T$-action to get a wonderful compactification $\Won(\T)$ of $\T$.  The unique closed $G \times G_\T$-orbit (where $G_\T = \Isom_G(\T,\T)$ is the $\T$-twisted inner form of $G$) is then $G/B \times \T/B$, where $\T/B$ parameterizes reductions of structure groups of $\T$ to a Borel.  Since $\T$ has a reduction of structure to a Borel, then $\T/B$ has a $K$-point.  Thus
the closed subscheme $G/B\times T/B$ has a $K$-point $s_0$. Now, using  Hensel's lemma, take a formal
deformation of this $K$-point of $\Won(T)$ to a $K[[x]]$-point $s$ whose generic fiber $s_\eta$ is in the interior $\T$ of $\Won(\T)$. Since the pullback of $\T$ to $\Spec K((x))$ has the rational point $s_\eta$, the pullback torsor is trivial. Thus, by Serre-Grothendieck conjecture over DVR \cite{nis}, the pullback of T is trivial over $\Spec K[[x]]$.  By restricting to the closed point $\Spec K$, the original torsor $\T$ is trivial.\qed

\proof[Proof of Corollary \ref{serre2}] Since $G$ is quasisplit, there exists a Borel subgroup $B$ defined over $k(S)$. For any $G$-torsor $E$, we define the twisted full flag $k(S)$-varieties $E/B$. The elementary obstruction of $E/B$ vanishes by \cite[Lem. 6.4]{Gille} and \cite[Lemma 2.2 (vi)]{BCTS}. Thus Theorem \ref{intromain} implies that the torsor $E$ admits a reduction to $B$. 

Let $Z$ be the center of $G$. Let $G'=G/Z$ be the adjoint form of $G$. For any $G$-torsor $\T$, by the first paragraph, the induced $G'$-torsor $\T'$ admits a reduction to $B'=B/Z$. By Lemma \ref{jasontrick}, $\T'$ is a trivial $G'$-torsor. Thus by long exact sequence of Galois cohomology, the torsor $\T$ admits a reduction to the center $Z$. 
\qed

\newcommand{\etalchar}[1]{$^{#1}$}


\begin{thebibliography}{ABD{\etalchar{+}}64}
	
	\bibitem[ABD{\etalchar{+}}64]{SGA3}
	M.~Artin, J.~E. Bertin, M.~Demazure, P.~Gabriel, A.~Grothendieck, M.~Raynaud,
	and J.-P. Serre.
	\newblock {\em Sch\'emas en groupes.}, volume 1963/64 of {\em S\'eminaire de
		G\'eom\'etrie Alg\'ebrique de l'Institut des Hautes \'Etudes Scientifiques}.
	\newblock Institut des Hautes \'Etudes Scientifiques, Paris, 1964.
	
	\bibitem[Art74]{artinversal}
	M.~Artin.
	\newblock Versal deformations and algebraic stacks.
	\newblock {\em Invent. Math.}, 27:165--189, 1974.
	
	\bibitem[BB73]{BBdecomp}
	A.~Bia{\l}ynicki-Birula.
	\newblock Some theorems on actions of algebraic groups.
	\newblock {\em Ann. of Math. (2)}, 98:480--497, 1973.
	
	\bibitem[BCTS08]{BCTS}
	M.~Borovoi, J.-L. Colliot-Th{\'e}l{\`e}ne, and A.~N. Skorobogatov.
	\newblock The elementary obstruction and homogeneous spaces.
	\newblock {\em Duke Math. J.}, 141(2):321--364, 2008.
	
	\bibitem[Beh]{behrendthesis}
	Kai Behrend.
	\newblock The lefschetz trace formula for the moduli stack of principal
	bundles.
	\newblock {\em PhD Thesis.}
	
	\bibitem[BGI71]{SGA6}
	P~Berthelot, A.~Grothendieck, and L.~editors. Illusie.
	\newblock {\em Th\'eorie des intersections et th\'eor\`eme de
		{R}iemann-{R}och}.
	\newblock Lecture Notes in Mathematics, Vol. 225. Springer-Verlag, Berlin,
	1971.
	\newblock S{\'e}minaire de G{\'e}om{\'e}trie Alg{\'e}brique du Bois-Marie
	1966--1967 (SGA 6), Dirig{\'e} par P. Berthelot, A. Grothendieck et L.
	Illusie. Avec la collaboration de D. Ferrand, J. P. Jouanolou, O. Jussila, S.
	Kleiman, M. Raynaud et J. P. Serre.
	
	\bibitem[BLR90]{BLR}
	Siegfried Bosch, Werner L{\"u}tkebohmert, and Michel Raynaud.
	\newblock {\em N\'eron models}, volume~21 of {\em Ergebnisse der Mathematik und
		ihrer Grenzgebiete (3) [Results in Mathematics and Related Areas (3)]}.
	\newblock Springer-Verlag, Berlin, 1990.
	
	\bibitem[Bor91]{Borel}
	Armand Borel.
	\newblock {\em Linear algebraic groups}, volume 126 of {\em Graduate Texts in
		Mathematics}.
	\newblock Springer-Verlag, New York, second edition, 1991.
	
	\bibitem[Bri10]{Brion}
	M~Brion.
	\newblock On automorphism groups of fiber bundles. preprint.
	\newblock 2010.
	
	\bibitem[Bri15]{Brion2015}
	Michel Brion.
	\newblock Which algebraic groups are {P}icard varieties?
	\newblock {\em Sci. China Math.}, 58(3):461--478, 2015.
	
	\bibitem[Car02]{Carrell1}
	James~B. Carrell.
	\newblock Torus actions and cohomology.
	\newblock In {\em Algebraic quotients. {T}orus actions and cohomology. {T}he
		adjoint representation and the adjoint action}, volume 131 of {\em
		Encyclopaedia Math. Sci.}, pages 83--158. Springer, Berlin, 2002.
	
	\bibitem[CTGP04]{CTGP}
	J.-L. Colliot-Th{\'e}l{\`e}ne, P.~Gille, and R.~Parimala.
	\newblock Arithmetic of linear algebraic groups over 2-dimensional geometric
	fields.
	\newblock {\em Duke Math. J.}, 121(2):285--341, 2004.
	
	\bibitem[CTS87]{CTS}
	Jean-Louis Colliot-Th{\'e}l{\`e}ne and Jean-Jacques Sansuc.
	\newblock La descente sur les vari\'et\'es rationnelles. {II}.
	\newblock {\em Duke Math. J.}, 54(2):375--492, 1987.
	
	\bibitem[DB81]{dubois}
	Philippe Du~Bois.
	\newblock Complexe de de {R}ham filtr\'e d'une vari\'et\'e singuli\`ere.
	\newblock {\em Bull. Soc. Math. France}, 109(1):41--81, 1981.
	
	\bibitem[dJHS11]{dJHS}
	A.~J. de~Jong, Xuhua He, and Jason~Michael Starr.
	\newblock Families of rationally simply connected varieties over surfaces and
	torsors for semisimple groups.
	\newblock {\em Publ. Math. Inst. Hautes \'Etudes Sci.}, (114):1--85, 2011.
	
	\bibitem[dJS03]{dJS1}
	A.~J. de~Jong and J.~Starr.
	\newblock Every rationally connected variety over the function field of a curve
	has a rational point.
	\newblock {\em Amer. J. Math.}, 125(3):567--580, 2003.
	
	\bibitem[dJS06]{dJS}
	A.~J. de~Jong and J.~Starr.
	\newblock Low degree complete intersections are rationally simply connected.
	preprint.
	\newblock 2006.
	
	\bibitem[FP97]{FP}
	W.~Fulton and R.~Pandharipande.
	\newblock Notes on stable maps and quantum cohomology.
	\newblock In {\em Algebraic geometry---{S}anta {C}ruz 1995}, volume~62 of {\em
		Proc. Sympos. Pure Math.}, pages 45--96. Amer. Math. Soc., Providence, RI,
	1997.
	
	\bibitem[GHS03]{GHS}
	Tom Graber, Joe Harris, and Jason Starr.
	\newblock Families of rationally connected varieties.
	\newblock {\em J. Amer. Math. Soc.}, 16(1):57--67 (electronic), 2003.
	
	\bibitem[Gil10]{Gille}
	Philippe Gille.
	\newblock Serre's conjecture {II}: a survey.
	\newblock In {\em Quadratic forms, linear algebraic groups, and cohomology},
	volume~18 of {\em Dev. Math.}, pages 41--56. Springer, New York, 2010.
	
	\bibitem[Gro62]{FGA}
	Alexander Grothendieck.
	\newblock {\em Fondements de la g\'eom\'etrie alg\'ebrique. [{E}xtraits du
		{S}\'eminaire {B}ourbaki, 1957--1962.]}.
	\newblock Secr\'etariat math\'ematique, Paris, 1962.
	
	\bibitem[Gro71]{SGA1}
	Alexander Grothendieck.
	\newblock {\em Rev\^etements \'etales et groupe fondamental}.
	\newblock Springer-Verlag, Berlin, 1971.
	\newblock S{\'e}minaire de G{\'e}om{\'e}trie Alg{\'e}brique du Bois Marie
	1960--1961 (SGA 1), Dirig{\'e} par Alexandre Grothendieck. Augment{\'e} de
	deux expos{\'e}s de M. Raynaud, Lecture Notes in Mathematics, Vol. 224.
	
	\bibitem[Gro05]{SGA2}
	Alexander Grothendieck.
	\newblock {\em Cohomologie locale des faisceaux coh\'erents et th\'eor\`emes de
		{L}efschetz locaux et globaux ({SGA} 2)}.
	\newblock Documents Math\'ematiques (Paris) [Mathematical Documents (Paris)],
	4. Soci\'et\'e Math\'ematique de France, Paris, 2005.
	\newblock S{\'e}minaire de G{\'e}om{\'e}trie Alg{\'e}brique du Bois Marie,
	1962, Augment{\'e} d'un expos{\'e} de Mich{\`e}le Raynaud. [With an
	expos{\'e} by Mich{\`e}le Raynaud], With a preface and edited by Yves Laszlo,
	Revised reprint of the 1968 French original.
	
	\bibitem[Har77]{Hartshorne}
	Robin Hartshorne.
	\newblock {\em Algebraic geometry}.
	\newblock Springer-Verlag, New Yor, 1977.
	\newblock Graduate Texts in Mathematics, No. 52.
	
	\bibitem[HS05]{HS}
	Joe Harris and Jason Starr.
	\newblock Rational curves on hypersurfaces of low degree. {II}.
	\newblock {\em Compos. Math.}, 141(1):35--92, 2005.
	
	\bibitem[Hum90]{Humphreys1}
	James~E. Humphreys.
	\newblock {\em Reflection groups and {C}oxeter groups}, volume~29 of {\em
		Cambridge Studies in Advanced Mathematics}.
	\newblock Cambridge University Press, Cambridge, 1990.
	
	\bibitem[Kem76]{kempfvanishing}
	George~R. Kempf.
	\newblock Linear systems on homogeneous spaces.
	\newblock {\em Ann. of Math. (2)}, 103(3):557--591, 1976.
	
	\bibitem[KK10]{lcdubois}
	J{\'a}nos Koll{\'a}r and S{\'a}ndor~J. Kov{\'a}cs.
	\newblock Log canonical singularities are {D}u {B}ois.
	\newblock {\em J. Amer. Math. Soc.}, 23(3):791--813, 2010.
	
	\bibitem[KM76]{detconstruction}
	Finn~Faye Knudsen and David Mumford.
	\newblock The projectivity of the moduli space of stable curves. {I}.
	{P}reliminaries on ``det'' and ``{D}iv''.
	\newblock {\em Math. Scand.}, 39(1):19--55, 1976.
	
	\bibitem[Kol86]{Kollarhigher1}
	J{\'a}nos Koll{\'a}r.
	\newblock Higher direct images of dualizing sheaves. {I}.
	\newblock {\em Ann. of Math. (2)}, 123(1):11--42, 1986.
	
	\bibitem[Kol96]{Kollar}
	J{\'a}nos Koll{\'a}r.
	\newblock {\em Rational curves on algebraic varieties}, volume~32 of {\em
		Ergebnisse der Mathematik und ihrer Grenzgebiete. 3. Folge. A Series of
		Modern Surveys in Mathematics [Results in Mathematics and Related Areas. 3rd
		Series. A Series of Modern Surveys in Mathematics]}.
	\newblock Springer-Verlag, Berlin, 1996.
	
	\bibitem[Kol03]{Kollar-Fund}
	J{\'a}nos Koll{\'a}r.
	\newblock Rationally connected varieties and fundamental groups.
	\newblock In {\em Higher dimensional varieties and rational points ({B}udapest,
		2001)}, volume~12 of {\em Bolyai Soc. Math. Stud.}, pages 69--92. Springer,
	Berlin, 2003.
	
	\bibitem[KP01]{KimRahul}
	B.~Kim and R.~Pandharipande.
	\newblock The connectedness of the moduli space of maps to homogeneous spaces.
	\newblock In {\em Symplectic geometry and mirror symmetry ({S}eoul, 2000)},
	pages 187--201. World Sci. Publ., River Edge, NJ, 2001.
	
	\bibitem[Lan52]{TsenLang}
	Serge Lang.
	\newblock On quasi algebraic closure.
	\newblock {\em Ann. of Math. (2)}, 55:373--390, 1952.
	
	\bibitem[Nis84]{nis}
	Yevsey~A. Nisnevich.
	\newblock Espaces homog\`enes principaux rationnellement triviaux et
	arithm\'etique des sch\'emas en groupes r\'eductifs sur les anneaux de
	{D}edekind.
	\newblock {\em C. R. Acad. Sci. Paris S\'er. I Math.}, 299(1):5--8, 1984.
	
	\bibitem[Pey04]{Peyre02}
	Emmanuel Peyre.
	\newblock Counting points on varieties using universal torsors.
	\newblock In {\em Arithmetic of higher-dimensional algebraic varieties ({P}alo
		{A}lto, {CA}, 2002)}, volume 226 of {\em Progr. Math.}, pages 61--81.
	Birkh\"auser Boston, Boston, MA, 2004.
	
	\bibitem[SdJ10]{dJS2}
	Jason Starr and Johan de~Jong.
	\newblock Almost proper {GIT}-stacks and discriminant avoidance.
	\newblock {\em Doc. Math.}, 15:957--972, 2010.
	
	\bibitem[Sko01]{SK}
	Alexei Skorobogatov.
	\newblock {\em Torsors and rational points}, volume 144 of {\em Cambridge
		Tracts in Mathematics}.
	\newblock Cambridge University Press, Cambridge, 2001.
	
	\bibitem[Sta10]{Starr}
	Jason~Michael Starr.
	\newblock Rational points of rationally simply connected varieties.
	\newblock In {\em Vari\'et\'es rationnellement connexes: aspects
		g\'eom\'etriques et arithm\'etiques}, volume~31 of {\em Panor. Synth\`eses},
	pages 155--221. Soc. Math. France, Paris, 2010.
	
	\bibitem[SX11]{starrxu}
	Jason Starr and Chenyang Xu.
	\newblock Rational points of rationally simply connected varieties over global
	function fields. preprint.
	\newblock 2011.
	
	\bibitem[Wit08]{Wittenberg0}
	Olivier Wittenberg.
	\newblock On {A}lbanese torsors and the elementary obstruction.
	\newblock {\em Math. Ann.}, 340(4):805--838, 2008.
	
\end{thebibliography}
\end{document}